\documentclass[11pt]{article}
\usepackage[margin=1in]{geometry}
\usepackage{amssymb,amsmath,enumerate, tikz, amsthm, amsfonts, mathtools}
\usepackage[T1]{fontenc}

\definecolor{burgundy}{rgb}{0.5, 0.0, 0.13}
\definecolor{darkgreen}{rgb}{0.0, 0.39, 0.0}
\definecolor{darkmagenta}{rgb}{0.55, 0.0, 0.55}
\definecolor{darkorange}{rgb}{1.0, 0.55, 0.0}
\definecolor{darkred}{rgb}{0.55, 0.0, 0.0}
\definecolor{darkgray}{rgb}{0.66, 0.66, 0.66}
\definecolor{lightblue}{rgb}{0.68, 0.85, 0.9}
\definecolor{lightgray}{rgb}{0.83, 0.83, 0.83}
\definecolor{mauve}{rgb}{0.88, 0.69, 1.0}
\definecolor{mauvetaupe}{rgb}{0.57, 0.37, 0.43}
\definecolor{navyblue}{rgb}{0.0, 0.0, 0.5}
\definecolor{palemagenta}{rgb}{0.98, 0.52, 0.9}

\usepackage[hidelinks, unicode, plainpages=false]{hyperref}
\usepackage{xcolor}
\hypersetup{
	colorlinks,
	linkcolor={burgundy},
	citecolor={mauvetaupe},
	urlcolor={navyblue}
}
\usepackage{cleveref}
\usepackage{todonotes}
\usepackage{thmtools}
\usepackage{setspace}
\newtheorem{theorem}{Theorem}[section]
\RenewCommandCopy{\theHtheorem}{\thetheorem}
\newtheorem{corollary}[theorem]{Corollary}
\RenewCommandCopy{\theHcorollary}{\thecorollary}
\newtheorem{proposition}[theorem]{Proposition}
\RenewCommandCopy{\theHproposition}{\theproposition}

\RenewCommandCopy{\theHconvention}{\theconvention}
\newtheorem{lemma}[theorem]{Lemma}
\RenewCommandCopy{\theHlemma}{\thelemma}
\newtheorem{definition}[theorem]{Definition}
\RenewCommandCopy{\theHdefinition}{\thedefinition}

\RenewCommandCopy{\theHclaim}{\theclaim}
\newtheorem{fact}[theorem]{Fact}
\RenewCommandCopy{\theHfact}{\thefact}

\RenewCommandCopy{\theHrequirement}{\therequirement}
\newtheorem{remark}[theorem]{Remark}
\RenewCommandCopy{\theHremark}{\theremark}

\RenewCommandCopy{\theHassumption}{\theassumption}
\newtheorem{observation}[theorem]{Observation}
\RenewCommandCopy{\theHobservation}{\theobservation}

\RenewCommandCopy{\theHnotation}{\thenotation}

\RenewCommandCopy{\theHconjecture}{\theconjecture}

\RenewCommandCopy{\theHexample}{\theexample}

\RenewCommandCopy{\theHconstruction}{\theconstruction}

\RenewCommandCopy{\theHconclusion}{\theconclusion}

\RenewCommandCopy{\theHquestion}{\thequestion}

\newcommand{\ift}{\text{if }}

\newcommand{\EE}{\mbox{\msbm E}}

\newcommand{\abs}[1]{\left| #1 \right|}
\newcommand{\comment}[1]{}

\newcommand{\indicate}[1]{\mathbf{1}\left\{ #1 \right\}}
\newcommand{\paren}[1]{\left( #1 \right)}
\newcommand{\set}[1]{\left\{ #1 \right\}}

\newcommand{\C}{\mathcal{C}}

\newcommand{\K}{\mathcal{K}}

\newcommand{\exVal}[1]{\mathbb{E}\left[ #1 \right]}
\newcommand{\ofS}[1]{\left[ #1 \right]}

\newcommand{\Hc}{\mathcal{H}}

\newcommand{\Oc}{\mathcal{O}}

\newcommand{\ve}{\varepsilon}

\newtheorem{prob}{Problem}
 
\DeclareMathOperator{\ver}{vert}
\DeclareMathOperator{\init}{init}
\newcommand{\maxvertdisc}{\delta_{\ver}^\infty}
\newcommand{\const}{\textcolor{black}{\alpha}}
\newcommand{\constt}{\textcolor{black}{\beta }}

\allowdisplaybreaks[4] 

\begin{document}

\date{November 5, 2026}

\title{An improved minimum degree threshold for fractional clique decompositions in partite graphs}

\author{Sean Longbrake\footnote{\textsc{Department of Mathematics, Emory University, Atlanta, United States.} The third author was partially supported by the National Science Foundation grant 2247013: Forbidden and Colored Subgraphs. \emph{Email addresses:} {\fontfamily{cmtt}\selectfont \{sean.longbrake, liana.yepremyan, xiangxiang.michael.zheng\}@emory.edu}}\footnotemark[1] \and  Michael Simkin\footnote{\textsc{Hebrew University of Jerusalem, Einstein Institute of Mathematics and Rachel and Selim Benin School of Computer Science and Engineering.} \emph{Email address:} {\fontfamily{cmtt}\selectfont michael.simkin@mail.huji.ac.il}} \and  Liana Yepremyan\footnotemark[1] \and Michael Zheng\footnotemark[1]}

\maketitle

\begin{abstract}
 We prove a Dirac-type condition for an  $r$-partite graph to admit a fractional clique decomposition: every sufficiently large $r$-partite graph with $n$ vertices in each part and partite minimum degree at least $(1 - \frac{1}{30r^2}) n$ admits a fractional $K_r$-decomposition, subject to necessary divisibility conditions. This implies the same upper bound on the asymptotic minimum degree threshold for $r$-partite graphs to admit a $K_r$-decomposition, using a well-known reduction.
\end{abstract}

\section{Introduction}

A $K_r$-decomposition of a graph is a partition of its edges into copies of $K_r$. A basic problem in extremal design theory is to determine how large the minimum degree must be in a graph to guarantee such a decomposition, subject to the necessary divisibility conditions. For triangles, Nash-Williams \cite{nash1970unsolved} conjectured that every sufficiently large $n$-vertex graph with minimum degree at least $3n/4$, even vertex degrees, and an edge count divisible by three admits a triangle decomposition. After much progress in the last half a century~\cite{gustavsson1991decompositions, yuster2005asymptotically, barber2016edge, garaschuk2014linear, dross2016fractional, dukes2020minimum, delcourt2021progress}, Delcourt and Postle \cite{delcourt2026proof} recently proved this
conjecture. A central ingredient is their resolution of its \emph{fractional relaxation}: every graph with minimum degree at least $3n/4$ admits a non-negative weighting of its triangles in which the total weight on each edge is one. 

More generally,  we say that a graph $G$ is $K_r$-\emph{divisible} if $e(G)$ is divisible by $\binom r2$ and every
vertex degree is divisible by $r-1$. The question then is to find the right minimum degree lower bound which guarantees that every $K_r$-divisible graph with that minimum degree admits a $K_r$-decomposition. Formally, let $\delta_{K_r}(n)$ be the least
$\delta>0$ such that every $K_r$-divisible graph $G$ with $n$ vertices
and $\delta(G)\geq\delta n$ has a $K_r$-decomposition.
The \emph{asymptotic $K_r$-decomposition threshold} is then defined as
$\displaystyle{
  \delta_{K_r}=\limsup_{n\to\infty}\delta_{K_r}(n).
}$

Next we define these parameters in the fractional setting. A \emph{fractional $K_r$-decomposition}  of $G$ is a function
$\omega\colon \mathcal K_r(G)\to\mathbb R_{\ge0}$ satisfying
\[
 \sum_{\substack{K\in\mathcal K_r(G): e\in E(K)}}\omega(K)=1
 \qquad\text{for every }e\in E(G),
\]

where $\K_r(G)$ denotes the set of copies of $K_r$ in $G$. Let $\delta_{K_r}^*(n)$ be the least $\delta>0$ such that
every $K_r$-divisible graph $G$ with $n$ vertices and $\delta(G)\geq\delta n$ has a
fractional $K_r$-decomposition. The \emph{asymptotic fractional
$K_r$-decomposition threshold} is
$\displaystyle{\delta_{K_r}^*=\limsup_{n\to\infty}\delta_{K_r}^*(n)}$. 

Fractional decompositions are useful both as objects in their own right and as a route to integral decompositions. For every fixed $r\geq3$, Glock, K\"uhn, Lo, Montgomery, and Osthus \cite[Corollary~1.2]{glock2019decomposition}
showed that the asymptotic minimum degree thresholds for fractional and integral $K_r$-decompositions coincide, that is, $\delta_{K_r}=\delta_{K_r}^*$. For $r\geq4$, Montgomery \cite{montgomery2019fractional} proved that  $\delta_{K_r}^* \leq 1-\frac{1}{100r}$. This bound on the asymptotic fractional
$K_r$-decomposition threshold is of the optimal order of magnitude in $r$. Indeed, for every $n$ divisible by $r(r-1)^2$, the balanced, complete $(r-1)$-partite
graph with $n/(r-1)$ vertices in each part is $K_r$-divisible, has minimum degree
$(1-\frac{1}{r-1})n$, and contains no $K_r$. The best lower bound construction \cite{gustavsson1991decompositions, yuster2005asymptotically} was actually recently improved. Indeed, a folklore generalization of Nash-Williams' conjecture, stated by Gustavsson in his 1991 thesis~\cite{gustavsson1991decompositions},  was recently proved to be false by Delcourt, Henderson, Lesgourgues, and Postle
\cite{delcourt2026counterexamples}. The conjecture said that  every
sufficiently large $K_r$-divisible graph $G$ with
$\delta(G)\ge(1-\frac{1}{r+1})n$ has a $K_r$-decomposition. The results in~\cite{delcourt2026counterexamples} disproved this conjecture in the stronger setting by showing that  for every $r\geq 4$, there exists some $c>1$ such that there exist infinitely many $K_r$-divisible graphs $G$ with minimum degree at least $\left(1-\frac{1}{c(r+1)}\right)n$ and no fractional $K_r$-decomposition. In fact, their results also imply the existence of such counterexamples even in the $r$-\emph{partite setting} which will be the focus of our paper.

This setting is closely related to
completion problems for combinatorial designs: triangle decompositions of
the complete balanced tripartite graph correspond to Latin squares, and
$K_r$-decompositions of the complete balanced $r$-partite graph correspond
to systems of $r-2$ mutually orthogonal Latin squares. Prescribing some
cliques and decomposing the remaining edges gives the corresponding
completion problem; see for example \cite{MR3663493, montgomery2017fractional, bowditch2019fractional}. Every $r$-clique in such an $r$-partite graph is \emph{transversal}, using
exactly one vertex from each part. The divisibility conditions in the $r$-partite setting are stronger, and to state them we need some definitions. 

Let $G$ be an $r$-partite graph with vertex partition
$V_1\cup\cdots\cup V_r$. We call $G$ \emph{balanced} if
$|V_1|=\cdots=|V_r|$, and define its \emph{partite minimum degree} by
\[
 \hat\delta(G)=\min_{i\ne j}\min_{v\in V_i} \deg_j(v)
\]
where $\deg_j(v) \coloneqq \abs{N(v) \cap V_j}$.\footnote{One could also consider our problem in a more general partite setting. Recently, Feng, Liu, and Yu \cite{feng2026fractional} studied sufficient partite minimum degree conditions for fractional $K_r$-decompositions of balanced $s$-partite graphs
with $s>r$.} In addition to equal part sizes, the partite setting has a local balancedness
condition that is necessary even for fractional decompositions. Indeed,
if $\omega$ is a fractional $K_r$-decomposition, then for every
$v\in V_i$ and every $j\ne i$,
\[
 \deg_j(v)
 =\sum_{u\in N(v)\cap V_j}
   \sum_{\substack{K\in\mathcal K_r(G)\\uv\in E(K)}}\omega(K)
 =\sum_{\substack{K\in\mathcal K_r(G)\\v\in V(K)}}\omega(K).
\]
Note that the final expression is independent of $j$, and thus $\deg_j(v)=\deg_k(v)$ for any $j,k\neq i$. Following Montgomery~\cite{montgomery2017fractional} and with slight abuse of wording, we say that $G$ is \emph{$K_r$-divisible} if
\begin{equation}\label{eq:intro-local-balance}
 \deg_j(v)=\deg_k(v)
 \quad\text{for all distinct }i,j,k\in[r]\text{ and every }v\in V_i.
\end{equation}
Throughout the partite setting, $K_r$-divisibility refers to
\eqref{eq:intro-local-balance}. It is easy to check that this condition also implies the usual
divisibility conditions. 

Similar to the non-partite setting, it was shown by Barber, Kühn, Lo, Osthus, and Taylor~\cite[Corollary 1.6]{MR3663493} that the asymptotic minimum degree thresholds for fractional and integral $K_r$-decompositions also coincide. We give the precise definition for these quantities below.

Fix $r\ge3$. Let $\hat{\delta}_{K_r}(n)$ be the least
$\delta>0$ such that every $K_r$-divisible $r$-partite graph $G$ with $n$ vertices in each part
and $\hat{\delta}(G)\geq\delta n$ admits a $K_r$-decomposition.
The \emph{asymptotic $r$-partite $K_r$-decomposition threshold} is then defined as
$\displaystyle{
  \hat{\delta}_{K_r}=\limsup_{n\to\infty}\hat{\delta}_{K_r}(n).
}$

Let $\hat\delta_{K_r}^*(n)$ be the least $\delta>0$ such that every balanced,
$K_r$-divisible, $r$-partite graph $G$ with $n$ vertices in each
part and $\hat\delta(G)\geq\delta n$ admits a fractional
$K_r$-decomposition. The \emph{asymptotic $r$-partite fractional $K_r$-decomposition threshold} is
\[
  \hat\delta_{K_r}^*
  =\limsup_{n\to\infty}\hat\delta_{K_r}^*(n).\]

As we stated earlier, it was shown in~\cite{MR3663493} that  $\hat\delta_{K_r}=\hat\delta_{K_r}^*$. Finding a $K_r$-decomposition in the $r$-partite setting turns out to be a harder problem than in the non-partite setting. Indeed, the following reduction was observed by Montgomery \cite{montgomery2017fractional}, and was also the passage to the disproof  of the folklore generalization of Nash-Williams' conjecture in the $r$-partite setting by~\cite{delcourt2026counterexamples}. Given a graph $G$, form $\widehat G$ from $r$ copies of $V(G)$,
joining vertices in distinct copies whenever their corresponding
vertices are adjacent in $G$. Then, $\widehat G$ is balanced, $K_r$-divisible, and satisfies
$\hat\delta(\widehat G)=\delta(G)$. Furthermore, $\widehat G$ admits a fractional
$K_r$-decomposition if and only if $G$ does. Consequently,
\[
  \hat\delta_{K_r}^*\ge\delta_{K_r}^*
  >1-\frac{1}{r+1}
  \qquad\text{for every }r\ge4.
\]

The best upper bound for general $r$ prior to our work was by Montgomery \cite{montgomery2017fractional} who showed that, for every integer $r\ge3$, $\hat\delta_{K_r}^*\le 1-\frac{1}{10^6r^3}$. Our main result improves Montgomery's upper bound to
$ \hat\delta_{K_r}^*\le 1-\frac{1}{30r^2}$, 
increasing the minimum degree deficit from order $r^{-3}$ to order $r^{-2}$ in the partite setting.

\begin{theorem}\label{thm:main}
For every integer $r\ge3$, there exists $n_0=n_0(r)$ such that the
following holds for all $n\ge n_0$. Every balanced, $K_r$-divisible,
$r$-partite graph $G$ with $n$ vertices in each part and
\[
 \hat\delta(G)\ge\left(1-\frac{1}{30r^2}\right)n
\]
admits a fractional $K_r$-decomposition.
\end{theorem}

For $r=3$, improving on work in \cite{daykin1983completion, chetwynd1985completing, gustavsson1991decompositions, bartlett2013completions, bowditch2019fractional, yu2026completion, he2026fractional}, Allsop, Bowtell, Lesgourgues, and Petrova
\cite{allsop2026improved} recently proved that
$\hat\delta_{K_3}^*<0.769$ which is better than our bound. Note that  our focus was on improving the dependency of  $\hat\delta_{K_r}^*$ on $r$, for large $r$. By the earlier mentioned result of Barber, K\"uhn, Lo, Osthus and Taylor~\cite{MR3663493}, our result implies the following corollary.

\begin{corollary}
    Let $r \geq 3$ be an integer and let $0<\ve < \frac{1}{30 r^{2}}$. There exists some $n_0 = n_0(\varepsilon,r)$ such that for every $n \geq n_0$, every balanced, $K_r$-divisible, $r$-partite graph $G$ with $rn$ vertices and $\hat\delta(G) \geq (1-\ve)n$ admits a $K_r$-decomposition.
\end{corollary}

\subsection{Organization}

\Cref{sec:sketch} discusses the general outline of our proof strategy. \Cref{sec:three} sets up some initial inequalities we will use. \Cref{sec:v-adjust} describes how we adjust the vertex weights to get a vertex-balanced weighting. \Cref{sec:e} then iteratively adjusts the edge weights. Importantly, the edge-adjusters do not change the vertex-balancedness. \Cref{sec:proof} proves that the limit of the weightings we constructed in \Cref{sec:e} assigns non-negative weights to all cliques.  

\section{Proof sketch}
\label{sec:sketch}
Let $G$ be a balanced, $r$-partite, $K_r$-divisible graph with parts $\abs{V_1} = \dots = \abs{V_r} = n$ satisfying $\hat{\delta}(G) \geq (1- \ve )n$. Recall that by $K_r$-divisibility, $\deg_j(v) = \deg_k(v)$ holds for distinct $i,j,k \in [r], v \in V_i$. In particular, the number of edges between any two parts will be the same: $e(V_i, V_j) = e(V_k, V_\ell)$ for any $i \neq j, k \neq \ell$. It will thus be convenient to define $\widehat{\deg(v)} = \deg_j (v)$ for $v \in V_i, j \neq i$, and $\widehat{e(G)} = e(G) / \binom{r}{2}$. Lastly, denote by $\mathcal{K}_r(G)$ the set of all $r$-cliques in $G$ and $k_r(G) = \abs{\mathcal K _r (G)}$.

In our proof, we will build a fractional $K_r$-decomposition for $G$ via an explicit limiting process. 

\begin{figure}[htb!]
    \centering
    \includegraphics[width=0.45\linewidth]{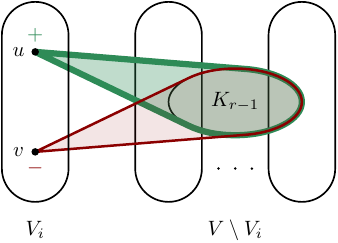}
    \caption{Vertex-adjuster increasing / decreasing the value of $u$ / $v$}
    \label{fig:vertex-adjuster}
\end{figure}

The starting point of this process will be the one in \cite{MR4651013}. We start by assigning each clique the weight $1$, so set $\omega_{\init}(K) = 1$ for all $r$-cliques $K$.
Note that this is in contrast to another potential starting point of assigning weight $\widehat{e(G)} / k_r(G)$ to all cliques. Under this normalization, we want to eventually find $\omega\colon \mathcal K _r(G) \to \mathbb R_{\geq 0 }$ which satisfies the following equation: 
\begin{equation*}
    \omega(e) \coloneqq \sum_{K\in \mathcal{K}_r(G):e\in E(K)}{\omega(K)}=\frac{ k_r(G)}{\widehat{e(G)}}
\end{equation*}
for every edge $e \in E(G)$. This is the condition for $\omega$ to be a fractional decomposition. Note, however, that it is necessary for $\omega$ to satisfy
\begin{equation}
    \omega(v) \coloneqq \sum_{K\in \mathcal{K}_r(G):v\in K}{\omega(K)} = \frac{k_r(G)}{\widehat{e(G)}} \cdot \widehat{\deg(v)} 
    \label{eq:vertex-good}
\end{equation} 
for any $v$ in $V(G)$. We will call such weightings \emph{vertex-balanced}. So in the first step, we aim towards the necessary condition by adjusting our initial weighting $\omega_{\init}$ to one that satisfies this property. For that, we will use \emph{vertex-adjusters} which only change the weights of two specified vertices in some part so that the weight is ``moved'' from one to the other, see \Cref{fig:vertex-adjuster}; in particular, the total amount of weight remains unchanged. Let us call the weighting that we find after vertex adjustment $\omega_0$.

Starting with $\omega_0$, we then build a sequence of vertex-balanced weightings $\omega_1, \omega_2, \dots $, where the difference for every edge $f$ between $\omega_i(f)$ and ${k_r(G)}/{\widehat{e(G)}}$ goes to zero with $i$ exponentially; we will call this difference the \emph{(edge-)discrepancy} of $f$ with respect to $\omega_i$ and denote it by $\delta_i(f)$. While we cannot guarantee that $\delta_i(f)$ always decreases by this rate for each edge $f$, we are able to guarantee that it is true for the \emph{maximum (edge-)discrepancy} 
\begin{equation}
    \delta_i^\infty \coloneqq \max_{e \in E(G)} \abs{\delta_i(e)}.
\end{equation}
Thus, if we were to take the limit of these weightings, we obtain a weighting $\omega_{\infty}$ of $G$ where the weight of every edge is exactly ${k_r(G)}/{\widehat{e(G)}}$. 

\begin{figure}[htb!]
    \centering
    \includegraphics[width=0.95\linewidth]{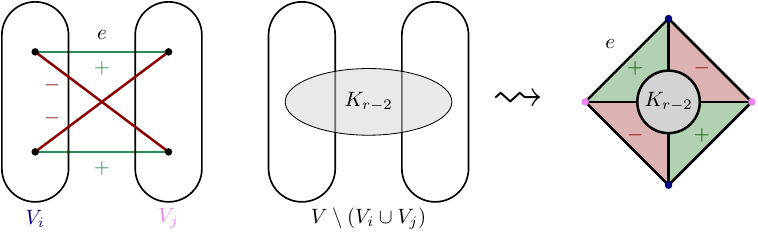}
    \caption{Edge-adjuster for an edge $e \in E(\textcolor{navyblue}{V_i}, \textcolor{palemagenta}{V_j})$}
    \label{fig:edge-adjuster}
\end{figure}

To build this sequence, we will make use of \emph{edge-adjusters}, whose supports are made of cycles of length four between two parts paired with a $K_{r - 2}$ such that any edge of the cycle induces a $K_r$ with the $K_{r-2}$, see \Cref{fig:edge-adjuster}. These edge-adjusters will preserve vertex-balancedness, while only adjusting weights on the edges of the cycle. To ensure that we only significantly adjust the weight of a specified edge $e$, we will take the sum of these edge-adjusters over all cycles $C$ containing $e$ and all $(r-2)$-cliques that are in the neighborhood of $C$. Using the fact that these weightings are vertex-balanced, we will be able to show that from $\omega_i$ to $\omega_{i+1}$, the weight of an edge $e$ moves towards $k_r(G) / \widehat{e(G)}$ by about $\delta_i(e) \pm \ve \delta_i^\infty$. This shows that the maximum discrepancy drops by about a factor of $\ve$, in particular that it goes to zero as $i \rightarrow \infty$. 

It then suffices to show that $\omega_i(K)$ is at least zero for every clique $K$ and index $i$. We do this by showing that the change in weight of $K$ for vertex balancing is no more than $\Oc(\ve r^2)$, and when passing from $\omega_i$ to $\omega_{i +1}$ for $i \geq 1$, the change of weight caused by the edge-adjusters is no more than $\Oc(\ve r^3 (c\ve)^i)$ for some absolute constant $c > 0$. 

However, 
applying the method naively for the transition from $\omega_0$ to $\omega_1$ as described above only allows for $\ve = \Oc(\frac{1}{r^3})$. 
For a vertex $v$ and an edge $uv$, let $k_r(v)$ and $k_r(uv)$ be the number of $r$-cliques that $v$ and $uv$ are contained in, respectively. For an edge $uv$, the discrepancy of this edge with respect to $\omega_0$ is roughly the edge's discrepancy with respect to $\omega_{\init}$, minus the discrepancies of its endpoints times $n^{-1}$
(see \Cref{genericvertexdiscrepencybound}); $$\delta_0(uv)\approx \underbrace{k_r(uv)-\frac{k_r(G)}{\widehat{e(G)}}}_{=\delta_{\init}(uv)} - \frac{\delta_{\ver} (u) + \delta_{\ver}(v)}{n} ,$$ where $\delta_{\ver}(u) = k_r(u)  - \frac{k_r(G)}{\widehat{e(G)}} \cdot \widehat{\deg(u)}$. One can show that 
\begin{equation*}
    \delta_{\ver} (u ) = \frac{1}{r-1} \sum_{v \in N(u)} \paren{k_r(uv) - \frac{k_r(G)}{\widehat{e
(G)}} };
\end{equation*}
note that $\exVal{k_r(e)} = k_r(G) / \widehat{e(G)}$ where $e \in E(G)$ is chosen uniformly at random. In particular, all of the above quantities involve bounding the difference $\abs{k_r(e) - k_r(f)}$ for edges $e,f \in E(G)$. Even in the case where $e,f$ are within the same parts, say $V_1$ and $V_2$, the best one can do in general is by writing the difference as a sum of indicators 
\begin{equation*}
    \sum_{v_3 \in V_3} \dots \sum_{v_r \in V_r} \paren{\indicate{e v_3 \dots v_r \in \mathcal K_r (G)} - \indicate{f v_3 \dots v_r \in \mathcal K_r (G)}}
\end{equation*}
and one-by-one splitting $V_i$ into $(V_i \cap (N(e) \cap N(f)))$ and $(V_i \cap (N(e)^c \cup N(f)^c))$. Using the fact that $\widehat{\deg(v)} \geq (1- \ve) n$ for any $v\in V(G)$, we obtain bounds of order $\mathcal O(\ve r n^{r-2})$, see again \Cref{genericvertexdiscrepencybound}. As the change of weight of a clique $K_r$ is roughly determined by the sum of the absolute values of the  discrepancies of its edges, we can only bound the change by $\Oc(\ve r^3)$ from $\omega_0$ to $\omega_1$ which does not give us the desired bound for $\ve$. 

To bypass this, we avoid cliques whose sum of absolute edge discrepancies is too large and exclude them from the first edge adjustment step $\omega_0 \rightsquigarrow \omega_1$ altogether. In other words, we will only consider edge-adjusters on $C_4$'s and $K_{r-2}$'s where the sum of edge discrepancies of the latter is bounded by $\Oc(\ve^2 r^4 n^{r-2})$. To get a reasonable bound on the number of such cliques, we initially clean $G$ so that every $r$-clique has a vertex of degree $(r - 1)(1 - \ve)n$, assuming without loss of generality that $(1 - \ve)n$ is an integer. This lets us assume that between any two parts, $G$ has no more than $(1 - \ve + \ve^2 r)n^2$ edges. After the cleaning, a simple probabilistic argument shows that this only excludes a small fraction of the edge-adjusters. As we only change the weight of each of these cliques by something roughly on the order of the discrepancies of the cycle times $ n^{- r}$, we obtain the desired bound.

\section{Setup}\label{sec:three}

We will now begin to explain in detail the procedure outlined in \Cref{sec:sketch}. In this section, we will establish all the remaining necessary definitions, define the preprocessing step and bound the initial discrepancies provided by our initial weighting $\omega_{\init}$. 

Throughout this section, let $G$ be a balanced, $r$-partite, $K_r$-divisible graph with $\hat{\delta}(G) = (1- \ve)n$ where $\ve \leq \const^{-1} r^{-2}$. \Cref{thm:main} is then equivalent to showing that $\const = 30$ is sufficient for $G$ to admit a fractional $K_r$-decomposition, provided $n$ is sufficiently large. It will be useful to extend $k_r(\cdot)$ from vertices and edges to arbitrary subsets of $V(G)$.
\begin{definition}
    For a set $S \subseteq V(G)$ let $\K_r(G, S)$ be the set of $K_r$'s in $G$ containing $S$. Furthermore, let $k_r(G, S) = \abs{\K_r(G, S)}$. We will suppress $G$ whenever the considered ambient graph is clear from context.
\end{definition}
Sometimes, we will need to consider $(r-t)$-cliques between $r-t$ specified parts.  
\begin{definition}
    For distinct $i_1, \dots, i_t \in [ r]$, let $G_{i_1, i_2, \dots, i_t} = G[V_{i} \colon i \in \set{i_1, \dots, i_t}]$ and $G^{i_1,i_2, \dots, i_t} = G[V_j \colon j \not \in \set{i_1, \dots, i_t}]$. Furthermore, let $k_{r-t}^{i_1, \dots, i_t}$ be the number of $(r-t)$-cliques in $G^{i_1, \dots, i_t}$ with the convention that $k_{0}^{i_1, \dots, i_r} = 1$. Lastly, for $Z \subseteq V(G^{i_1, \dots, i_t})$, let $k_{r-t}^{i_1, \dots, i_t}(Z)$ denote the number of $(r-t)$-cliques containing $Z$.
\end{definition}
Similarly, we need the common neighborhood of a set of vertices to a specified part.
\begin{definition}
    For $i \in [r]$ and $S \subseteq V(G)$, let $N_i(S) \coloneqq V_i \cap \bigcap_{s \in S} N(s)$ denote the \emph{$i$-neighborhood of $S$}; note that ${N_i(S)} = \emptyset$ if $S \cap V_i \neq \emptyset$. Furthermore, let $N_i^c(S) \coloneqq V_i \setminus N_i(S)$. 
\end{definition}

By our minimum degree condition, we note the following fact.  
\begin{fact}\label{eq:big}
    $G$ has $k_r(G) \geq n^r \cdot \prod_{i = 1}^{r -1} (1 - i \ve) \geq \paren{1 - r (r-1) \varepsilon /2} n^r$ many $K_r$'s. In particular, we may assume $k_r(G) \geq (1- (2\const)^{-1}) n^r$ provided that $\varepsilon \leq \const^{-1} r^{-2}$. 
\end{fact}

\subsection{Preprocessing}

In this section, we will reduce finding fractional $K_r$-decompositions in general $K_r$-divisible graphs to those that are additionally \emph{$(\ve, r)$-minimal}.
\begin{definition}
    A balanced, $r$-partite, $K_r$-divisible graph $G$ is \emph{$(\ve, r)$-minimal} if $\hat{\delta}(G) = (1-\ve)n$ and for every $K \in \K_r(G)$,
    \begin{equation*}
        \min_{v \in V(K)} \widehat{\deg(v)} = (1- \ve) n.
    \end{equation*}
\end{definition}

The next lemma justifies why it will be sufficient to only consider $(\ve, r)$-minimal $G$ for our algorithm where from now on, we will additionally assume that $\ve = \const^{-1} r^{-2}$. 
\begin{lemma}
    Let $\ve \in (0,1)$. Given a balanced, $r$-partite, $K_r$-divisible graph $G$ with $\hat{\delta}(G) \geq ( 1- \ve) n$, there exists a subgraph $G' \subseteq G$ that is $(\ve, r)$-minimal and if $G'$ has a fractional $K_r$-decomposition, then $G$ does as well.
\end{lemma}
\begin{proof}
    While there are cliques $K$ with $\min_{v \in V(K)} \widehat{\deg(v)}  > (1-\ve)n$, we may remove all the edges of $K$. Since the deletion of a clique reduces the degree of a vertex into a part only by one, we may iterate until the resulting graph $G'$ is $(\ve, r)$-minimal. Note that this preserves $K_r$-divisibility. If the resulting graph has a fractional $K_r$-decomposition, we may extend it to one of $G$ by simply assigning each $r$-clique in the removal weight one; all other removed $r$-cliques are assigned weight zero. By construction, this is a fractional $K_r$-decomposition.
\end{proof}

What makes minimality crucial is that it gives us an upper bound on the number of edges to pair with the lower bound given by the minimum degree.
\begin{proposition}\label{eq:global_edge}
    If $G$ is $(\ve,r)$-minimal, then 
    \begin{equation*}
        (1- \ve)n^2 \leq \widehat{e(G)} \leq ( 1- \ve 
         + \ve^2 (r-1))n^2. 
    \end{equation*}
\end{proposition} 
    
\begin{proof}
    It suffices to prove the upper bound; the lower bound follows directly from the handshaking lemma. Suppose that for every $j \in [r]$, the number of $v \in V_j$ satisfying $\widehat{\deg(v)} > (1 - \ve)n$ is larger than $\ve (r-1) n$. Then, there is a $K_{r}$ satisfying that every vertex has partite degree strictly greater than $(1 - \ve)n$; one may greedily obtain such a $K_r$. Thus, there is some $j$ such that the set $X_j\coloneqq  \{v \in V_j: \widehat{\deg(v)} > (1 - \ve)n\}$ satisfies $|X_j| \leq \ve (r-1) n$. It follows then for all $i \neq j$ that $\sum_{u \in V_j}|\deg_i(u) - (1 - \ve)n| \leq |X_j| \ve n$, and so 
    \begin{equation*}
        \widehat{e(G)} = e(G_{ij}) \leq (1 - \ve + \ve^2 (r-1))n^2. \qedhere
    \end{equation*}
\end{proof}
This estimate will be crucial in proving \Cref{lem:avg-clique-disc}.

\subsection{Bounding initial discrepancies}

Before getting into vertex adjustment in \Cref{sec:v-adjust}, we first make some observations about how far away the constant weighting $\omega_{\init}$ is from being a fractional $K_r$-decomposition. To make this formal, we will consider the edge-discrepancies of $\omega_{\init}$, denoted $\delta_{\init}$, and the vertex-discrepancy $\delta_{\ver}$. 
In particular, we establish bounds on the \emph{maximum edge-discrepancy (with respect to $\omega_{\init}$)} $\delta_{\init}^\infty \coloneqq \max_{e\in E(G)} \abs{\delta_{\init}(e)}$ and the \emph{maximum vertex-discrepancy} $\delta_{\ver}^\infty \coloneqq \max_{v \in V(G)} \abs{\delta_{\ver} (v)}$. 
Recall from \Cref{sec:sketch} that for $v \in V(G), e \in E(G)$,
\begin{align*}
    \delta_{\ver} (v) &= \omega_{\init} (v) - \frac{k_r(G)}{\widehat{e(G)}} \cdot \widehat{\deg(v)}, & \delta_{\init}(e) &= \omega_{\init} (e) - \frac{k_r(G)}{\widehat{e(G)}}.
\end{align*}

Lastly, throughout this section, note that $\omega_{\init} (e) = k_r(e)$ and $\omega_{\init}(v) = k_r(v)$.

\begin{lemma}
    For $e,e'\in E(G_{ij})$, $1 \leq i < j \leq r$, we get 
    \begin{align*}
        \abs{k_r(e) - k_r(e')} &\leq 2 (r-2) \ve n ^{r-2}.
    \end{align*}
    In particular, $\abs{\delta_{\init}(e)} = \abs{k_r(e) - k_r(G) / \smash{\widehat{e(G)}}}\leq 2 (r-2) \ve n^{r-2}$. 
    \label{clm:e-diff}
\end{lemma}

\begin{proof}
    Without loss of generality assume that $(i,j) = (1,2)$ and $k_r(e) \geq k_r(e')$. We may expand as follows:
    \begin{align*}
        {}& k_r(e) - k_r(e') \\
        ={}& \sum_{v_3 \in V_3} \dots \sum_{v_r \in V_r} \paren{\indicate{ev_3 \dots v_r \in \K_r(G)} - \indicate{e'v_3 \dots v_r \in \K_r(G)}} \\
        ={}&  {\sum_{v_3 \in N_3(e) \cap N_3(e')} \dots \sum_{v_r \in V_r}  \paren{\indicate{e v_3 \dots v_r \in \K_r(G)} - \indicate{e'v_3 \dots v_r \in \K_r(G)}}} \\
        & + \underbrace{\sum_{v_3 \in ( N_3(e) \triangle N_3(e'))} \dots \sum_{v_r \in V_r}  \paren{\indicate{e v_3\dots v_r \in \K_r(G)} - \indicate{e' v_3 \dots v_r \in \K_r(G)}}}_{\leq 2\varepsilon n^{r-2}} \\ 
        \vdotswithin{=} & \\ 
        \leq{}&  \sum_{v_3 \in N_3(e) \cap N_3(e')} \dots \sum_{v_r \in  N_r(e) \cap N_r(e')}  \paren{\indicate{e v_3 \dots v_r \in \K_r(G)} - \indicate{e'v_3 \dots v_r \in \K_r(G)}}  \\
        & + 2\varepsilon (r-2) n^{r-2} \\
        ={}& 2\varepsilon (r-2)n^{r-2}.
    \end{align*}
    The inequality in the second equation holds since the difference of indicators is positive only if $v_3 \in N_3(e) \setminus N_3(e')$; $\abs{N_3(e) \setminus N_3(e')}$ is at most $2\ve n$ and upper bounding for such $v_3$ the indicator by one gives us $n^{r-3}$ per $v_3$. The last equality follows from the fact that if $v_i \in V_i \cap N(e) \cap N(e')$ for $3 \leq i \leq r$, then $e v_3 \dots v_r$ forms an $r$-clique if and only if $v_3\dots v_r$ forms an $(r-2)$-clique. This happens if and only if $e' v_3 \dots v_r$ forms an $r$-clique.

    Finally, note that $\exVal{k_r(f)} = k_r(G) / \widehat{e(G)}$ if $f \in E(G_{ij})$ is chosen uniformly at random, so
    \begin{equation*}
        {\delta_{\init}(e)} = {k_r(e)-  \exVal{k_r(f)}} = \exVal{k_r(e) - k_r(f)} \in [- 2\ve (r-2) n^{r-2}, 2 \ve (r-2) n^{r-2}]
    \end{equation*}
    by linearity of expectation and since the latter bound holds for any arbitrary $f\in E(G_{ij})$.
\end{proof}

For some parts of our proof, more refined estimates on $k_r(e)$ will be necessary, see \Cref{lem:avg-clique-disc}. For this, we will use an analogue of \cite[Lemma 4.2 (ii)]{barber2017fractional}, see \Cref{scaryinequality}.
Roughly speaking, this lemma states that every $(r-2)$-clique that doesn't extend an edge $xy$ to an $r$-clique must contain a vertex outside the neighborhood of $x$ or $y$, and each such vertex extends in about $k_{r - 3}^{1, 2, i}$ many ways to an $(r-2)$-clique.

Similar to \Cref{clm:e-diff}, we first show 
\begin{lemma}
    \label{clm:sub-clq-diff}
    For $1 \leq i_1, i_2, \dots, i_{t+1} \leq r$ distinct, and $z \in V_{i_{t+1}}$, we have 
    \begin{equation*}
        k_{r- t -1}^{i_1,\dots, i_{t+1}} - k_{r - t}^{i_1, \dots, i_t} (z) \leq \ve (r- t - 1) n^{r-t - 1}.
    \end{equation*}
\end{lemma}
\begin{proof}
    Without loss of generality, assume that $\set{i_1, \dots, i_t} = [r-t+1,r]$ and $i_{t+1} = r-t$. Clearly,
    \begin{align*}
        & k_{r- t -1}^{i_1,\dots, i_{t+1}} - k_{r - t}^{i_1, \dots, i_t} (z) \\
        ={}& \sum_{v_1 \in V_{1}} \sum_{v_2 \in V_2} \dots \sum_{{ v_{r-t-1} \in V_{r-t-1}}} \paren{\indicate{v_1\dots v_{r-t-1}\in \K_{r-t-1}(G)}  - \indicate{z v_1\dots v_{r-t-1}\in \K_{r-t}(G)} } \\
        ={}& \sum_{v_1 \in N_1(z)} \sum_{v_2 \in V_2} \dots \sum_{{ v_{r-t-1} \in V_{r-t-1}}} \paren{\indicate{v_1\dots v_{r-t-1}\in \K_{r-t-1}(G)}  - \indicate{z v_1\dots v_{r-t-1}\in \K_{r-t}(G)} } \\
        &+  \underbrace{\sum_{v_1 \in V_1 \setminus N_1(z)}  \dots \sum_{{ v_{r-t-1} \in V_{r-t-1}}} \paren{\indicate{v_1\dots v_{r-t-1}\in \K_{r-t-1}(G)}  - \indicate{z v_1\dots v_{r-t-1}\in \K_{r-t}(G)} } }_{\leq \ve n^{r-t-1}}\\ 
        \vdotswithin{=}& \\
        \leq{}& \sum_{v_1 \in N_1(z)} \dots \sum_{{ v_{r-t-1} \in N_{r-t-1}(z)}} \paren{\indicate{v_1\dots v_{r-t-1}\in \K_{r-t-1}(G)}  - \indicate{z v_1\dots v_{r-t-1}\in \K_{r-t}(G)} } \\
        & + \ve (r-t-1) n^{r-t-1} \\
        ={}& \ve (r-t-1) n^{r-t-1}.
    \end{align*}
    As $v_i \in N_i(z)$ for all $i \in [r-t-1]$, the last equality follows since $v_1 \dots v_{r-t-1} \in \K_{r-t-1}(G)$ if and only if $z v_1 \dots v_{r-t-1} \in \K_{r-t}(G)$.
\end{proof}

We now use \Cref{clm:sub-clq-diff} to prove our desired inequality. 

\begin{lemma}\label{scaryinequality}
Suppose $xy \in E(G_{i_1i_2}), 1 \leq i_1 , i_2 \leq r$ distinct. Then,
 $$ 0 \leq {k_r(xy) - k_{r - 2}^{i_1, i_2} + \sum_{j \in [r] \setminus \set{i_1, i_2}}|N^c_j(x) \cup N^c_j(y)|k_{r -3}^{i_1,i_2,j}} \leq 4(\ve (r-2))^2 n^{r - 2}.$$
\end{lemma}
\begin{proof}
    Without loss of generality suppose $i_1 = 1, i_2 = 2$.
    Using the inclusion-exclusion principle, we see that 
    \begin{equation*}
        k_{r - 2}^{1, 2} - k_r(xy) = \sum_{i = 1}^{r-2} (-1)^{i-1} \sum_{Z \subseteq (N^c(x) \cup N^c(y)) \setminus (V_1 \cup V_2), |Z| = i}k_{r - 2}^{1, 2}(Z).
    \end{equation*}
    Moving the $i =1$-term and applying the Bonferroni inequality \cite{bonferroni1936teoria} yields that 
    \begin{align*}
         0 \leq {k_r(xy) - k_{r - 2}^{1, 2}  + \sum_{z \in( N^c(x) \cup N^c(y)) \setminus (V_1 \cup V_2)}k_{r - 2}^{1, 2}(z) }&\leq  {\sum_{Z \subseteq (N^c(x) \cup N^c(y)) \setminus (V_1 \cup V_2), |Z| = 2}k_{r - 2}^{1, 2}(Z)}\\
        &\leq \frac{(r-2)^2}{2} (2 \ve n)^2  n^{r - 4} \\
        &= 2 (\ve (r-2))^2 n^{r-2}.
    \end{align*}

    Using \Cref{clm:sub-clq-diff}, we thus get
    \begin{align*}
        0 \leq{}& {k_r(xy) - k_{r - 2}^{1, 2}  + \sum_{i \in [3, r]}|N^c_i(x) \cup N^c_i(y)|k_{r -3}^{1,2,i} }  \\
        \leq{}& {k_r(xy) - k_{r - 2}^{1, 2} + \sum_{z \in( N^c(x) \cup N^c(y)) \setminus (V_1 \cup V_2)}k_{r - 2}^{1, 2}(z) + \sum_{i \in [3, r]}\sum_{z \in N^c_i(x) \cup N^c_i(y)}\left( k_{r - 3}^{1, 2, i} - k_{r - 2}^{1, 2}(z)\right) }  \\
        \leq{}& 2 (\ve (r-2))^2n^{r - 2} + {\sum_{i \in [3, r]}\sum_{z \in N^c_i(x) \cup N^c_i(y)}\left( k_{r - 3}^{1, 2, i} - k_{r - 2}^{1, 2}(z)\right)}\\
        \leq{}& 4 (\ve (r-2))^2n^{r - 2} 
    \end{align*}
    as desired.
\end{proof}

For $v \in V(G)$, recall that the vertex-discrepancy of $v$ is defined by $$\delta_{\ver}(v) = k_r(v) - \frac{k_r(G)}{\widehat{e(G)}} \cdot \widehat{\deg(v)}.$$

\begin{lemma}
    For $v \in V(G)$, we have 
    \begin{equation*}
        \abs{\delta_{\ver}(v)} 
        \leq 2 \ve (r-2) n^{r-1}.
    \end{equation*}
    In particular, $\maxvertdisc \leq 2 \ve (r-2) n^{r-1}$.
    \label{lem:maxvertdisc}
\end{lemma}
\begin{proof}
    Without loss of generality assume $v \in V_1$. Note that
    \begin{align*}
        \delta_{\ver}(v) &= \paren{\sum_{u\in N_2(v)} k_r(uv)} - \frac{k_r(G)}{\widehat{e(G)}} \cdot \widehat{\deg(v)} \\
        &= \sum_{u \in N_2(v)} \paren{k_r(uv) - \frac{k_r(G)}{\widehat{e(G)}}} \\
        &= \sum_{u \in N_2(v)} \delta_{\init}(uv).
    \end{align*}
    Hence, we obtain the claim using the triangle inequality and \Cref{clm:e-diff}.
\end{proof}
\begin{remark}\label{rem:v-disc_e-disc}
    In general, one can show that for $v \in V_i, i \in [r],$ we have 
    \begin{align*}
        \delta_{\ver}(v) = \sum_{u \in N_j(v)} \paren{k_r(uv) - \frac{k_r(G)}{\widehat{e(G)}}} = \sum_{u \in N_j(v)} \delta_{\init}(uv)
    \end{align*}
    for any $j \in [r] \setminus \set{i}$.
\end{remark}

Similar to \Cref{clm:e-diff}, one can show that the initial weights of any two vertices do not differ by much. 
\begin{lemma}    \label{clm:v-diff}
    For $i \in [r], u, v \in V_i$, $\abs{k_r(v) - k_r(u)} \leq \varepsilon (r-1) n^{r-1}$. In particular, 
    \begin{equation*}
        \abs{\delta_{\ver}(v) - \delta_{\ver}(u)} \leq \ve r n^{r-1}.
    \end{equation*}
\end{lemma}
\begin{proof}
    Without loss of generality we may relabel such that $i = 1$ and $k_r(v) \geq k_r(u)$.

    We may expand as follows:
    \begin{align*}
        {}& k_r(v) - k_r(u) \\
        ={}& \sum_{v_2 \in V_2} \sum_{v_3 \in V_3} \dots \sum_{v_r \in V_r} \paren{\indicate{vv_2 \dots v_r \in \K_r(G)} - \indicate{uv_2 \dots v_r \in \K_r(G)}} \\
        ={}& \sum_{v_2 \in V_2 \cap N(u,v)} {\sum_{v_3 \in V_3} \dots \sum_{v_r \in V_r}  \paren{\indicate{vv_2 \dots v_r \in \K_r(G)} - \indicate{uv_2 \dots v_r \in \K_r(G)}}} \\
        & + \underbrace{\sum_{v_2 \in V_2 \cap( N(u) \triangle N(v))} \sum_{v_3 \in V_3} \dots \sum_{v_r \in V_r}  \paren{\indicate{vv_2 \dots v_r \in \K_r(G)} - \indicate{uv_2 \dots v_r \in \K_r(G)}}}_{\leq \varepsilon n^{r-1}} \\ 
        \vdotswithin{=} & \\ 
        \leq{}&  \sum_{v_2 \in V_2 \cap N(u, v)} \dots \sum_{v_r \in V_r \cap N(u,v)}  \paren{\indicate{vv_2 \dots v_r \in \K_r(G)} - \indicate{uv_2 \dots v_r \in \K_r(G)}}  \\
        & + \varepsilon (r-1) n^{r-1} \\
        ={}& \varepsilon (r-1)n^{r-1}.
    \end{align*}
    The last equality follows from the fact that if $v_i \in V_i \cap N(u,v)$ then $u v_2 \dots v_r$ forms an $r$-clique if and only if $v_2\dots v_r$ forms an $(r-1)$-clique. This happens if and only if $v v_2 \dots v_r$ forms an $r$-clique. Lastly, as $(1-\ve) n \leq \widehat{\deg(u)}, \widehat{\deg(v)} \leq n$ and $k_r(G) \leq \sum_{e \in E(G_{12})} k_r(e) \leq \widehat{e(G)} n^{r-2}$, the second claim follows.
\end{proof}

Finally, as $G$ is $K_r$-divisible, we get the following simple, yet effective
\begin{lemma}
    For any $i \in [r]$, $\sum_{x_i \in V_i} \delta_{\ver}(x_i) = 0$.
    \label{obs:vert-disc-zero}
\end{lemma}
\begin{proof}
    From $G$ being $K_r$-divisible, we immediately get
    \begin{align*}
        \sum_{v_i \in V_i} \delta_{\ver}(v_i) = \sum_{v_i \in V_i} k_r(v_i) - \frac{ k_r(G)}{\widehat{e(G)}} \underbrace{\sum_{v_i \in V_i} \widehat{\deg(v_i)}}_{= \widehat{e(G)}} 
        = 0
    \end{align*}
    as desired.
\end{proof}

\section{Phase 1: Vertex adjustment}
\label{sec:v-adjust}

Throughout this section, let $G$ be a balanced, $r$-partite, $K_r$-divisible graph that is $(\ve, r)$-minimal where $\ve = \const^{-1} r^{-2}$. In this section, we will define the preliminary weighting $\omega_0$ on the $K_r$'s of $G$. Unlike the initial weighting $\omega_{\init}$, $\omega_0$ will be vertex-balanced. We show that we can achieve this while not deviating too much from $\omega_{\init}$. In particular, $\omega_0$ assigns positive weights provided that $\ve$ is small.

\begin{definition}
    Let $x_i, y_i \in V_i$ be a pair of distinct vertices. For $R \in \K_r(G)$, let $$\omega_{x_i, y_i}(R) \coloneqq \frac{1}{h(x_i, y_i)} \sum_{H \in \Hc(x_i, y_i)}(\indicate{R = H_{x_i}} - \indicate{R = H_{y_i}})$$where $\mathcal{H}(x_i, y_i)$ is the set of ``overlapping'' $K_r$'s $H = (H_{x_i}, H_{y_i}) \in \K_r(x_i) \times \K_r(y_i)$ such that $\abs{V(H_{x_i}) \cap V(H_{y_i})} = r-1$, and $h(x_i, y_i) \coloneqq \abs{\Hc(x_i, y_i)}$. 
\end{definition}
Note that as $\hat{\delta}(G)$ is large, we have for any pair $x_i, y_i$ many such overlapping $r$-cliques.
\begin{fact}\label{h-big} 
    Given distinct $x_i, y_i \in V_i$, $i \in [r]$, $h(x_i, y_i) \geq \paren{1 - r (r+1) \varepsilon/2} n^{r-1}$. In particular, $h(x_i, y_i) \geq (1 - \const^{-1}) {n^{r-1}}$ provided that $\ve \leq \const^{-1} r^{-2}$.
\end{fact}
From $\hat{\delta}(G)$ being large, it actually follows that varying one argument does not change $h$ by much.
\begin{lemma}
    For distinct $x, u,v \in V_i, i \in [r],$ we have $\abs{h(x,u) - h(x,v)} \leq \ve (r-1)n^{r-1}$.
    \label{clm:h-diff}
\end{lemma}
We use a proof similar to \Cref{clm:e-diff}.
\begin{proof}
    Without loss of generality we may relabel such that $i = 1$ and $h(x,v) \geq h(x,u)$. As $(H_x, H_v) \in \Hc(x,v)$ if and only if $V(H_v) \setminus \set{v} \subseteq N(x)$, $H_v \in \K_r(v)$ and $V(H_x) = V(H_v) \; \triangle \; \set{x,v}$, we may write 
    \begin{align*}
        {}& h(x,v) - h(x,u) \\
        ={}& \sum_{v_2 \in N_2(x) } \dots \sum_{v_r \in N_r(x) } \paren{\indicate{vv_2 \dots v_r \in \K_r(G)} - \indicate{uv_2 \dots v_r \in \K_r(G)}} \\
        ={}& \sum_{v_2 \in N_2(x, u,v)} {\dots \sum_{v_r \in N_r(x)}  \paren{\indicate{vv_2 \dots v_r \in \K_r(G)} - \indicate{uv_2 \dots v_r \in \K_r(G)}}} \\
        & + \underbrace{\sum_{v_2 \in N_2(x) \cap( N(u) \triangle N(v))} \dots \sum_{v_r \in N_r(x)}  \paren{\indicate{vv_2 \dots v_r \in \K_r(G)} - \indicate{uv_2 \dots v_r \in \K_r(G)}}}_{\leq \varepsilon n^{r-1}} \\ 
        \vdotswithin{=} & \\ 
        \leq{}&  \sum_{v_2 \in N_2(x, u, v)} \dots \sum_{v_r \in N_r(x, u,v)}  \paren{\indicate{vv_2 \dots v_r \in \K_r(G)} - \indicate{uv_2 \dots v_r \in \K_r(G)}}  \\
        & + \varepsilon (r-1) n^{r-1} \\
        ={}& \varepsilon (r-1)n^{r-1}.
    \end{align*}
    The last equality follows from the fact that if $v_i \in  N_i(x,u,v)$ then $u v_2 \dots v_r$ forms an $r$-clique if and only if $v_2\dots v_r$ forms an $(r-1)$-clique. This happens if and only if $v v_2 \dots v_r$ forms an $r$-clique.
\end{proof}

The idea behind $\omega_{x_i, y_i}$ is to shift the weight from, say, $y_i$ to $x_i$ if $x_i$ has a negative deficiency and $y_i$ a positive deficiency. While this can be done using a single $H \in \mathcal H(x_i, y_i)$, to not change the weight too much we average over all such $H$. Overall, we see the following
\begin{observation}
    $\omega_{x_i, y_i}(v) = \begin{cases} 1 & \ift  v = x_i\\
- 1 & \ift   v = y_i \\
0 & \text{otherwise}.
\end{cases}$
\label{obs:overlapping}
\end{observation}

Of course, to establish vertex-balancedness, it is not guaranteed that  pairing vertices with opposite discrepancies would be possible. In particular, one does not expect the existence of an $x_i, y_i$ with $\delta_{\ver}(x_i) + \delta_{\ver}(y_i) = 0$. On the other hand, we do know that $\sum_{x_i \in V_i} \delta_{\ver}(x_i) = 0$ by \Cref{obs:vert-disc-zero}. Thus, we make use of these $\omega_{x_i, y_i}$ over all pairs $x_i, y_i \in V_i$. 

The contribution of each gadget  will be determined for all $i \in [r]$ by 
\begin{definition}
    For $i \in [r]$ and $R \in \K_r(G)$, let
    \begin{equation*}
    \omega_{V_i}(R) \coloneqq - \frac{1}{2n}\sum_{(x_i, y_i) \in V_i^2\colon x_i \neq y_i} (\delta_{\ver}(x_i) - \delta_{\ver}(y_i)) \cdot \omega_{x_i, y_i}(R). 
\end{equation*}

Furthermore, set
\begin{align*}
    \omega_0 \coloneqq \omega_{\init} + \sum_{i \in [r]} \omega_{V_i}.
\end{align*}

\end{definition}
We first show that $\omega_0$ is indeed vertex-balanced and that it does not decrease too much in weight for any $r$-clique.
\begin{lemma}
    For $v \in V(G)$, $\omega_0(v) = \frac{k_r(G)}{ \widehat{e(G)}} \cdot \widehat{\deg(v)}$.
    \label{obs:omegaVi}
\end{lemma}
\begin{proof}
    It suffices to show that 
    \begin{equation*}
        \omega_{V_i} (v) = \begin{cases}
        -\delta_{\ver}(v), & v \in V_i \\
        0, &\text{otherwise}
    \end{cases}
    \end{equation*}
    for $v \in V(G), i \in [r]$. Indeed, this would imply that 
    \begin{align*}
        \omega_0(v) = k_r(v) - \delta_{\ver}(v) =\frac{k_r(G)}{\widehat{e(G)}} \cdot \widehat{\deg(v)}.
    \end{align*}
    
    By \Cref{obs:overlapping}, we have that $\omega_{V_i}(v) = 0$ for $v \not \in V_i$; for $v \in V_i$, \Cref{obs:overlapping} gives
    \begin{align*}
        \omega_{V_i} (v) &=  - \frac{1}{2n} \sum_{u \in V_i \setminus \set{v}} \ofS{(\delta_{\ver}(u) - \delta_{\ver}(v)) \cdot \omega_{u,v}(v) + (\delta_{\ver}(v) - \delta_{\ver}(u)) \cdot \omega_{v,u}(v)} \\
        &= \frac{1}{n} \sum_{u \in V_i \setminus \set{v}} \paren{\delta_{\ver}(u) - \delta_{\ver}(v)} \\ 
        &= - \delta_{\ver}(v)
    \end{align*}
    where the last equality is by \Cref{obs:vert-disc-zero}.
\end{proof}

\begin{lemma}
    For any $r$-clique $R \in \K_r(G)$, $\omega_0(R) \geq 1 -  (\const - 1)^{-1}$ provided that $\const > 1$.
    \label{lem:omega_0_positive}
\end{lemma}

\begin{proof}
    Let $R \in \K_r(G)$ be an $r$-clique with vertices $z_i \in V_i \cap V(R), i \in [r]$. We calculate
\begin{align*}
    \omega_0(R) &= 1 + \sum_{i \in [r]}\omega_{V_i}(R) \\
     &= 1 - \frac{1}{n} \sum_{i \in [r]}\sum_{y_i \in V_i \setminus \set{z_i}} \frac{(\delta_{\ver}(z_i) - \delta_{\ver}(y_i))}{h(z_i, y_i)} \indicate{y_i \in N(V(R) \setminus \{z_i\})} \\
     &\geq 1 - \frac{\ve r^2}{1 - \const^{-1}}  \tag{\Cref{clm:v-diff}, \Cref{h-big}}\\
     &\geq 1 - \frac{1}{\const - 1 } 
\end{align*}
where the latter inequality holds for $\ve \leq \const^{-1} r^{-2}$.
\end{proof}
We end this section by giving bounds on the edge-discrepancy of $\omega_0$, see \Cref{genericvertexdiscrepencybound}. We first show that going from $\omega_{\init}$ to $\omega_0$, the change in weight is determined by the vertex-discrepancy of the edge's endpoints. By the triangle inequality, it follows that the edge-discrepancy of an edge $e$ with respect to $\omega_0$ is bounded by  the initial edge-discrepancy of $e$, plus the vertex-discrepancy of its endpoints times $n^{-1}$, and an error term of order $\mathcal{O}(\ve^2r^2n^{r-2})$. 
\begin{lemma}\label{boundvertex0}
    For any $e = x_1x_2 \in E(G)$,
    \begin{equation*}
        |\omega_0(e) - k_r(e)| \leq \frac{1}{(1 - \const^{-1}) n} \abs{ \delta_{\ver}(x_1)} + \frac{1}{(1 - \const^{-1}) n} \abs{\delta_{\ver}(x_2)} + \frac{4 \ve^2 r (r-1)n^{r-2} }{(1-\const^{-1})^2} .
    \end{equation*}
\end{lemma}
\begin{proof}
    Without loss of generality suppose $x_i \in V_i$ for $i \in [2]$. Fix $v_1^\ast\in V_1 \setminus \set{x_1}$ and $v_2^* \in V_2 \setminus \set{x_2}$. As 
\begin{align*}
    \omega_0(e) &= k_r(e) + \sum_{i \in [r]}\underbrace{\omega_{V_i}(e)}_{= 0 \text{ for }i \not \in [2]}\\
    &= k_r(e) - \frac{1}{2n} \sum_{i \in \set{1,2}} \sum_{(u_i, v_i) \in V_i^2 \colon u_i \neq v_i} \paren{\delta_{\ver} (u_i) - \delta_{\ver}(v_i)} \cdot \underbrace{\omega_{u_i, v_i}(e)}_{= 0 \text{ if } x_i \not \in \set{u_i, v_i}}\\
    &= k_r(e) - \frac{1}{n} \sum_{i \in \set{1,2}} \sum_{v_i \in V_i \setminus \set{x_i}} \paren{\delta_{\ver} (x_i) - \delta_{\ver}(v_i)} \cdot \omega_{x_i, v_i}(e)\\
    &= k_r(e) - \frac{1}{n} \sum_{R \in \K_r(e)}\sum_{i \in \set{1, 2}}\sum_{v_i \in V_i \setminus \set{x_i}} \frac{\delta_{\ver}(x_i) - \delta_{\ver}(v_i)}{h(x_i, v_i)}\sum_{H \in \Hc(x_i, v_i)}\indicate{R = H_{x_i}} \\
    &= k_r(e) - \frac{1}{n} \sum_{R \in \K_r(e) }\sum_{i \in \set{1, 2}}\sum_{v_i \in V_i \setminus \set{x_i}} \frac{\delta_{\ver}(x_i) - \delta_{\ver}(v_i)}{h(x_i, v_i)} \indicate{V(R)\; \triangle \; \set{x_i, v_i} \in \mathcal{K}_r(G)}, 
\end{align*}
we see that
\begin{align*}
    |\omega_0(e) - k_r(e)|
    ={}& \abs{ \frac{1}{n} \sum_{R \in \K_r(e) }\sum_{i \in \set{1, 2}}\sum_{v_i \in V_i \setminus \set{x_i}} \frac{\delta_{\ver}(x_i) - \delta_{\ver}(v_i)}{h(x_i, v_i)} \indicate{V(R) \; \triangle \; \set{x_i, v_i} \in \mathcal{K}_r(G)}}\\
    \leq{}& \abs{ \frac{1}{n} \sum_{R \in \K_r(e) }\sum_{i \in \set{1, 2}}\sum_{v_i \in V_i \setminus \set{x_i}} \frac{\delta_{\ver}(x_i) - \delta_{\ver}(v_i)}{h(x_i, v_i)}} \\ 
    &+ \abs{ \frac{1}{n} \sum_{R \in \K_r(e) }\sum_{i \in \set{1, 2}}\sum_{v_i \in V_i \setminus \set{x_i}} \frac{\delta_{\ver}(x_i) - \delta_{\ver}(v_i)}{h(x_i, v_i)} \indicate{v_i \not \in N(V(R) \setminus \set{x_i} )}} \\
    \leq{}&  \abs{ \frac{1}{n} \sum_{R \in \K_r(e) }\sum_{i \in \set{1, 2}}\sum_{v_i \in V_i \setminus \set{x_i}} \frac{\delta_{\ver}(x_i) - \delta_{\ver}(v_i)}{h(x_i, v_i)}} + \frac{2\ve^2 r (r-1) n^{r-2}}{(1 - \const^{-1})} \tag{$\abs{V_i \cap N^c(V(R) \setminus \set{x_i})} \leq  \ve (r-1) n$, \Cref{clm:v-diff}, \Cref{h-big}}\\
    \leq{}& \abs{ \frac{1}{n} \sum_{R \in \K_r(e) }\sum_{i \in \set{1, 2}}\sum_{v_i \in V_i} \frac{\delta_{\ver}(x_i) - \delta_{\ver}(v_i)}{h(x_i, v_i^*)}} +\frac{2 \ve^2 r(r-1) n^{r-2}}{(1- \const^{-1})} \\ &+  \abs{ \frac{1}{n}\sum_{R \in \K_r(e) } \sum_{i \in \set{1,2}} \sum_{v_i \in V_i \setminus \set{x_i}}\frac{\delta_{\ver}(x_i) - \delta_{\ver}(v_i)}{h(x_i, v_i)h(x_i, v_i^*)}\left(h(x_i, v_i) - h(x_i, v_i^*) \right)} \\
    \leq{}& \abs{ k_r(e) \left( \frac{\delta_{\ver}(x_1)}{h(x_1, v_1^*)} + \frac{\delta_{\ver}(x_2)}{h(x_2, v_2^*)} \right)} +  \frac{4 \ve^2 r (r-1)n^{r-2} }{(1-\const^{-1})^2}  \tag{\Cref{clm:v-diff}, \Cref{obs:vert-disc-zero}, \Cref{h-big}, \Cref{clm:h-diff}}\\
    \leq{}& \frac{1}{(1 - \const^{-1}) n} \paren{\abs{ \delta_{\ver}(x_1)} + \abs{\delta_{\ver}(x_2)}} + \frac{4 \ve^2 r (r-1)n^{r-2} }{(1-\const^{-1})^2} 
\end{align*}
which was to be shown.
\end{proof}

\begin{corollary}\label{genericvertexdiscrepencybound} For every $uv \in E(G)$,
\begin{align*}
    |\delta_0(uv)| &\leq \abs{\delta_{\init}(uv)} +  \frac{1}{(1 - \const^{-1}) n} \abs{ \delta_{\ver}(u)} + \frac{1}{(1 - \const^{-1}) n} \abs{\delta_{\ver}(v)} + \frac{4 \ve^2 r (r-1)n^{r-2} }{(1-\const^{-1})^2}. 
\end{align*} 
    In particular, $\delta_0^\infty \leq 13 \ve (r-2) n^{r-2} / 2$ provided $\const \geq 15$. 
\end{corollary}
\begin{proof}
    By definition, $\delta_0(uv) = \omega_0(uv) - \frac{k_r(G)}{\widehat{e(G)}}$. 
    By \Cref{boundvertex0}, 
    \begin{align*}
    |\delta_0(uv)| &= \abs{\omega_0(uv) - \frac{k_r(G)}{\widehat{e(G)}}}\\
    &\leq \abs{k_r(uv) - \frac{k_r(G)}{\widehat{e(G)}}} + |\omega_0(uv) - k_r(uv)|\\
    &\leq \abs{\delta_{\init}(uv)} +  \frac{1}{(1- \const^{-1})n} \abs{ \delta_{\ver}(u)} + \frac{1}{(1- \const^{-1})n} \abs{\delta_{\ver}(v)} + \frac{4 \ve^2 r (r-1)n^{r-2} }{(1-\const^{-1})^2} .
    \end{align*}
    
    For the bound on $\delta_0^\infty$, apply \Cref{clm:e-diff} on $\abs{\delta_{\init}(uv)}$ and \Cref{lem:maxvertdisc} for $\abs{\delta_{\ver}(u)}, \abs{\delta_{\ver}(v)}$:
    \begin{align*}
        \delta_0^\infty &\leq 2 \ve (r-2) n^{r-2} + \frac{4 \ve (r-2) n^{r-2}}{1 -  \const^{-1}} + \frac{4 \ve^2 r (r-1)n^{r-2} }{(1-\const^{-1})^2}  \\
        &\leq \paren{1 + \frac{2}{1 - \const^{-1} }  + \frac{2r \ve }{(1- \const^{-1})^2} \frac{r-1}{r-2}} 2 \ve (r-2) n^{r-2} \\
        &= \paren{1 + \frac{2}{1 - \const^{-1} }  + \frac{2}{\const r (1- \const^{-1})^2} \frac{r-1}{r-2} } 2 \ve (r-2) n^{r-2} \\
        &\leq \frac{13}{2} \ve (r-2)n^{r-2}. \qedhere{}
    \end{align*}
\end{proof}

\renewcommand{\EE}{\mathbb{E}}

\section{Phase 2: Edge adjustment}\label{sec:e}
Let $G$ be a balanced, $r$-partite, $K_r$-divisible graph that is $(\ve, r)$-minimal where $\ve = \const^{-1} r^{-2}$. In this section, we will define $\omega_{i+1}$ given the vertex-balanced weighting $\omega_i$, $i \geq 0$. Note that $\omega_0$ is defined in \Cref{sec:v-adjust}. We will show that this new weighting $\omega_{i+1}$ is also vertex-balanced and has maximum edge-discrepancy $\delta_{i+1}^\infty = \Oc(\ve \delta_i^\infty)$ where 
\begin{align*}
    \delta_{i}^\infty &=  \displaystyle{\max_{e \in E(G)} \abs{ \delta_i(e)}}, & \delta_i(e) &= \omega_{i} (e) - \frac{k_r(G)}{\widehat{e(G)}},
\end{align*}
provided that $\varepsilon \ll 1$, see also \Cref{lem:conv}. Our final candidate will then be $\omega_\infty = \displaystyle{\lim_{i \to \infty} \omega_i}$.

For the proof that $\omega_{i}(K)$ does not decrease too much in $i$ for $K \in \K_r(G)$, we extend $\delta_i$ to arbitrary subsets of edges.
\begin{definition}
     For $S \subseteq E(G)$, set
\begin{equation*}
    \delta_i(S) \coloneqq \sum_{e \in S} \abs{\delta_i (e)}, 
\end{equation*}
where we use the convention that $\delta_i(H) = \delta_i(E(H))$ for $H \subseteq G$. 
\end{definition}
\subsection{The definition of edge-adjusters and \texorpdfstring{$\omega_{i+1}$}{the next weighting}}

For the edge-adjusters, it is crucial to \emph{not} use $r$-cliques $R$ where all the edges need a lot of adjustments. For example, in the worst case scenario it may be that $\delta_0 (R) \sim \binom{r}{2} \delta_0^\infty \lesssim \ve r^3 n^{r-2}$ by \Cref{genericvertexdiscrepencybound}. This would result in a bound where $\ve = C/r^3$ as outlined in \Cref{sec:sketch}.

We will show in \Cref{lem:avg-clique-disc} that this worst case scenario does not occur on average; we defer its proof to \Cref{sec:technical_lemma}.

\begin{lemma}\label{lem:avg-clique-disc}
    Let $xy \in E(G)$. If $\const \geq 26$, then $$\frac{1}{k_r(xy)} \sum_{R \in  \K_r(xy)} {\delta_0(R - x-y) }\leq 15\ve^2 (r-1)^2 (r-2)^2 n^{r-2}.$$
\end{lemma}

Since we want edge-adjusters to preserve vertex-balancedness while changing the weight of some specified edge $e \in E(G_{ij})$ (and not too many others), the basic idea will be to consider $4$-cycles $C \subseteq G_{ij}$ containing $e$ and $K_{r-2}$'s on the remaining parts that extend each edge of $C$ to a clique. Assigning $r$-cliques formed along the cycle weights $\pm 1$ in an alternating pattern preserves vertex-balancedness and only leads to changes in weight for $C$'s edges. Naturally, instead of just considering one such cycle $C$, we will average the contribution over many such cycles containing $e$ and $K_{r-2}$'s to keep the contribution of $e$'s adjusters to other edges small. Moreover, as alluded to in the beginning, we will exclusively consider $(r-2)$-cliques that extend to useful cliques only. In the following, $\constt $ will be  a properly chosen constant greater than one. 

\begin{definition}
    We  call an $(r-2)$-clique $W$ \emph{useful} if $\delta_0(W) \leq \constt  \cdot 15 \ve^{2} (r-1)^2 (r-2)^{2} n^{r - 2}$.
\end{definition}

\begin{definition}
    We call an $r$-clique $R$ \emph{useful with respect to an edge $e$} if $R=W\cup \{e\}$ for some useful $(r-2)$-clique $W$. We denote the collection of useful cliques by $\widetilde{\K_r(e)}$.
\end{definition}

\begin{remark}
    Using~\Cref{genericvertexdiscrepencybound}, every useful $r$-clique  satisfies:
    \begin{align*}
        \delta_0 (R) &\leq \constt \cdot 15 \varepsilon^{2} (r-1)^2 (r-2)^{2} n^{r-2} + 2(r-1)\paren{ \frac{13}{2} \ve (r-2) n^{r-2}} \\
        &\leq \paren{15\constt \ve (r-1) (r-2) + 13} \ve (r-1)(r-2)n^{r-2} \\ 
        &\leq \paren{\frac{15 \constt}{\const} + 13}\ve (r-1)(r-2) n^{r-2}.
    \end{align*}
\end{remark}

\begin{definition}
    For a $4$-cycle $C \subseteq G_{ij}, i,j \in [r]$,  let ${N(C)}$ be the $(r-2)$-cliques $W = \{v_\ell \in V_\ell:  \ell \neq i, j\}$ such that $W$ is useful, $W \cup \set{e}$ induces an $r$-clique for all $e \in E(C)$.
\end{definition}

Note that by \Cref{lem:avg-clique-disc} and Markov's inequality we know that 
\begin{corollary}    \label{cor:useful_W}
    The number of useful $r$-cliques containing an edge $xy \in E(G)$ is at least 
    \begin{equation*}
        \abs{\widetilde{\K_r(xy)}} \geq \paren{1 - \frac{1}{\constt}}k_r(xy) \geq \paren{1- \frac{1}{2\const}} \paren{1 - \frac{1}{\constt}}n^{r-2}.
    \end{equation*}
    
    In particular for any $4$-cycle $C \subseteq G_{ij}$, we have 
    \begin{equation*}
        \abs{{N(C)}} \geq \paren{1 - \frac{1}{\const} - \frac{2}{\constt}} n^{r-2}.
    \end{equation*}
\end{corollary}
\begin{proof}
    Apply Markov's inequality on $\mathbb{P}(X > \constt \cdot \exVal{X})$ where $X = \delta_0(R - x - y)$ and $R \in \K_r(xy)$ is chosen uniformly at random. The first inequality then follows from \Cref{lem:avg-clique-disc}. 

    For the second inequality, note that for an $(r-2)$-clique $W$ to form an $r$-clique with each $e \in E(C)$, it suffices that $W$ induces an $r$-clique with a matching $e_1, e_2 \in E(C)$. Hence, by inclusion-exclusion,
    \begin{equation*}
        \abs{N(C)} \geq \abs{\widetilde{\K_r(e_1)}} + \abs{\widetilde{\K_r(e_2)}} - n^{r-2} \geq \paren{1 - \frac{1}{\const} - \frac{2}{\constt}} n^{r-2} 
    \end{equation*}
    as desired.
\end{proof}

\begin{definition}
    For $e \in E(G_{ij})$, let $\C(e)$ be the set of $4$-cycles $C \subseteq G_{ij}$ containing $e$.\footnote{Similarly, for $S \subseteq V_i \cup V_j \cup E(G_{ij})$, we let $\C(S)$ be the $4$-cycles $C \subseteq G_{ij}$ containing $S$.} For $e \in E(G)$ and $C \in \C(e)$, let  
\begin{equation*}
    \psi_{C, e}(R) \coloneqq \frac{1}{\abs{N(C)}} \sum_{W\in N(C)} \sum_{f \in E(C)} (-1)^{\sigma(e,f)} \cdot \indicate{R = G[W \cup f]}, 
\end{equation*}
where $\sigma(e,f) \coloneqq \abs{e \cap f} \pmod{2}$ denotes the parity of $f$ with respect to $e$. Furthermore, let 
\begin{align*}
    \psi_e  \coloneqq \frac{1}{\abs{\C(e)}} \sum_{C \in \C(e)} \psi_{C,e}.
\end{align*}
\end{definition}

It is easy to see 
\begin{fact}\label{fact:v-balanced-preserved}
    $\psi_e(v) = 0$ for every vertex $v \in V(G)$ and $e \in E(G)$.
\end{fact}
In particular, adjusting using $\psi_e$ leaves a vertex-balanced function vertex-balanced.

\begin{fact}
    \label{fact:sum_C_f}
    Let $e \in E(G)$ and $C \in \C(e)$. Then we have 
    \begin{align*}
        \psi_{C, e}(f) = \begin{cases}
            (-1)^{\sigma(e,f)}, & f \in E(C) \\
            0, & \text{otherwise}.
        \end{cases}
    \end{align*}
\end{fact}

\begin{definition}
    Given the vertex-balanced weighting $\omega_i$, we (formally) set 
\begin{align*}
    \omega_{i+1} \coloneqq{}& \omega_i - \sum_{e \in E(G)} \delta_{i}(e) \psi_{e} \\
    ={}& \omega_i -  \sum_{e \in E(G)}  \frac{\delta_{i}(e)}{\abs{\C(e)}} \sum_{C \in \C(e)} \psi_{C,e}, \\
    \omega_{\infty} \coloneqq{}& \lim_{i \to \infty} \omega_i.
\end{align*}
\end{definition}

While the definition of $\omega_{i+1}$ is more technical, we inductively get from \Cref{fact:v-balanced-preserved} the following
\begin{corollary}\label{cor:seq:v-balanced}
    $\omega_{i}$ is vertex-balanced for all $i \geq 0$.
\end{corollary}

We wish to show that, after normalization, $\omega_{\infty}$ is a desired fractional $K_r$-decomposition. To do so, after proving \Cref{lem:avg-clique-disc} in \Cref{sec:technical_lemma}, we first show in \Cref{lem:conv} that the maximum edge-discrepancy tends to zero at the desired rate which in particular implies that $\omega_{\infty}$ is well-defined; note that this part of the argument goes through regardless of whether the sequence of weightings actually stays non-negative. We then finish this section by providing estimates on how much the weight of a clique decreases from $\omega_i$ to $\omega_{i+1}$.

For both sections, we will make use of the following representation of $\omega_{i+1}(R)$.
\begin{lemma}\label{lem:omega_i+1_simplified}
    For $i \geq 0$ and $R \in \K_r(G)$, we have 
    \begin{equation*}
        \omega_{i+1}(R) = \omega_i(R) - \sum_{f \in E(R)}\sum_{C \in \C(f)}\sum_{e \in E(C)} 
       \frac{(-1)^{\sigma(e,f)}\delta_i(e)\indicate{\exists W \in N(C): G[W \cup f] = R}}{|N(C)||\C(e)|}.
    \end{equation*}
\end{lemma}

\begin{proof}
    A direct calculation shows that 
\begin{align*}
    \omega_{i+ 1}(R) 
    ={}& \omega_i(R) - \sum_{e \in E(G)} \delta_i(e) \psi_e(R)\\
    ={}& \omega_i(R) - \sum_{e \in E(G) } \frac{\delta_i(e)}{\abs{\C(e)}} \sum_{C \in \C(e)} \psi_{C,e}(R) \\
    ={}& \omega_i(R) - \sum_{e \in E(G) } \frac{\delta_i(e)}{\abs{\C(e)} } \sum_{C \in \C(e)}\sum_{W\in N(C)} \frac{1}{\abs{N(C)}}\sum_{f \in E(C)} (-1)^{\sigma(e,f)} \indicate{R = G[W \cup f]}\\
    ={}& \omega_i(R) - \sum_{f \in E(R)}\sum_{C \in \C(f)}\sum_{e \in E(C)} 
       \frac{(-1)^{\sigma(e,f)}\delta_i(e)\indicate{\exists W \in N(C): G[W \cup f] = R}}{|N(C)||\C(e)|}
\end{align*}
where the last equality follows from the fact that $\indicate{R = G[W \cup f]} = 1$ only if $f \in E(R)$.
\end{proof}

\subsection{Proof of \texorpdfstring{\Cref{lem:avg-clique-disc}}{technical lemma}}

\label{sec:technical_lemma}
In this subsection, we will use \Cref{scaryinequality} to prove the aforementioned lemma which shows that the average edge-discrepancy of the clique is $\Oc(\ve^2r^4 n^{ r - 2})$.

\begin{proof}[Proof of \Cref{lem:avg-clique-disc}]

Without loss of generality we may assume that $x \in V_1, y\in V_2$. Note that for $r = 3$, the bound is immediate, so we will from now on assume $r \geq 4$. Double counting gives
\begin{align*}
    \frac{1}{k_r(xy)} \sum_{R \in \K_r(xy)} \delta_0(R -x-y) \tag{\Cref{eq:big}}
    ={}& \sum_{uv \in E(G^{1,2})} |\delta_0(uv)| \underbrace{\frac{k_r(uvxy)}{k_r(xy)}}_{\leq \frac{n^{r-4}}{\paren{1 - \const^{-1}/2} n^{r-2}}}  \\ 
    \leq{}& \frac{1}{\paren{1 - \frac{1 }{2 \const}}n^2} \sum_{uv \in E(G^{1,2})} \abs{\delta_0(uv)}. \tag{$\ast$} \label{eq:hard_lem_1}
\end{align*}

Recall from \Cref{genericvertexdiscrepencybound} that 
\begin{equation*}
    \abs{\delta_0(uv)} \leq \abs{\delta_{\init}(uv)}+ \frac{1}{1- \const^{-1}} \ofS{ \frac{\abs{\delta_{\ver}(u)}}{n}  +  \frac{\abs{\delta_{\ver} (v)} }{n}+ \frac{4}{1 - \const^{-1}}\ve^2 r(r-1) n^{r-2}}. 
\end{equation*}
Hence, as
\begin{align*}
    \frac{1}{n} \sum_{uv \in E(G^{1,2})} \paren{\abs{\delta_{\ver} (u)} +\abs{ \delta_{\ver}(v)}} &= \frac{1}{n} \sum_{u \in V(G^{1,2})} \sum_{v \in N(u) \setminus (V_1 \cup V_2)} \abs{\delta_{\ver} (u)} \\
    &\leq \sum_{u \in V(G^{1,2})} \sum_{v \in N(u) \setminus (V_1 \cup V_2)}  \frac{\frac{1}{r-3}\sum_{w \in N(u) \setminus (V_1 \cup V_2)} \abs{\delta_{\init}(uw)}}{n} \\
    &\leq \sum_{u \in V(G^{1,2})} \sum_{w \in N(u) \setminus (V_1 \cup V_2)} \abs{\delta_{\init}(uw)} \\
    &= 2 \sum_{uv \in E(G^{1,2})} \abs{\delta_{\init}(uv)} \\
\end{align*}

by \Cref{rem:v-disc_e-disc}, we may upper bound the sum $\sum_{uv \in E(G^{1,2})} \abs{\delta_0(uv)}$ in \eqref{eq:hard_lem_1} by 
\begin{align*}
    & \frac{3}{1 - \const^{-1}} \sum_{uv \in E(G^{1,2})} \abs{\delta_{\init}(uv)} + \frac{4\binom{r-2}{2}}{(1- \const^{-1})^2 } \ve^2 r (r-1)  n^{r} \\
    \leq{}& \frac{3}{1 - \const^{-1}} \sum_{ij \in \paren{[r] \setminus \set{1,2}}^{(2)}}\sum_{uv \in E(G_{ij})} \abs{k_r(uv) - \frac{1}{\widehat{e(G)}}\sum_{ab \in E(G_{ij})}k_{r}(ab)}  + \frac{2  \ve^2 r(r-2)^3 n^{r} }{(1- \const^{-1})^2}  \\ 
    \leq{}& \frac{3}{(1 - \const^{-1}) \widehat{e(G)}} \sum_{ij \in \paren{[r] \setminus \set{1,2}}^{(2)}} \sum_{ab, uv \in E(G_{ij})}  \abs{k_r(uv) - k_{r}(ab)} + \frac{2  \ve^2 r (r-2)^3 n^{r} }{(1- \const^{-1})^2}. 
    \tag{$\ast\ast$} \label{eq:hard_lem_2}
\end{align*}

    By \Cref{scaryinequality}, the following holds: 
    \begin{align*}
        & \abs{k_r(uv) - k_r(ab)} \\
        \leq{}& \abs{ \sum_{k \neq i, j} k_{r - 3}^{i,j,k}(|N^c_k(u)\cup N^c_k(v)| - |N^c_k(a) \cup N^c_k(b)|)}  + 4(\ve (r-2))^2 n^{r - 2} \\
        \leq{}& \sum_{k \neq i, j} k_{r - 3}^{i, j,k}  \big( \abs{ |N^c_k(u)\cup N^c_k(v)| - |N_k^c(u)| - |N^c_k(v)| + |N^c_k(a)| + |N_k^c(b)| - |N^c_k(a) \cup N^c_k(b)| }\\
        & \ \ \ \ \ \ \ \ \ \ \ \ \ + \abs{|N_k^c(u)| + |N^c_k(v)| - |N^c_k(a)| - |N_k^c(b)|} \big)  + 4(\ve (r-2))^2 n^{r - 2} \\
        ={}& \sum_{k \neq i, j} k_{r - 3}^{i, j,k} \big( \abs{|N^c_k(u)\cap N^c_k(v)| - |N^c_k(a) \cap N^c_k(b)| } + \abs{|N_k^c(u)| + |N^c_k(v)| - |N^c_k(a)| - |N_k^c(b)|} \big) \\
        & \ \ \ \ \ \ \ \ \ \ \ \ \  + 4(\ve (r-2))^2 n^{r - 2}. \tag{$\ast\ast\ast$} \label{eq:hard_lem_3}
    \end{align*}   

    Before plugging this into \eqref{eq:hard_lem_2}, let us first deal with the terms in the large parentheses individually: For the first part, we obtain
    \begin{align*}
        &  \sum_{ab, uv \in E(G_{ij})} \sum_{k \neq i, j} k^{i,j,k}_{r-3} \abs{|N^c_k(u)\cap N^c_k(v)| - |N^c_k(a) \cap N^c_k(b)|} \\
        \leq{}& n^{r-3} \sum_{ab, uv \in E(G_{ij})} \sum_{k \neq i, j} \paren{|N^c_k(u)\cap N^c_k(v)| + |N^c_k(a) \cap N^c_k(b)|} \\
        \leq{}& 2n^{r-1} \sum_{ab \in E(G_{ij})} \sum_{k \neq i, j}|N^c_k(a) \cap N^c_k(b)| \\
        \leq{}& 2n^{r-1} \sum_{k \neq i,j} \sum_{a \in V_i} \sum_{w \in N^c_k(a)} \abs{N^c_j (w)} \\
        \leq{}& 2 \ve^2 (r-2) n^{r+2}. \tag{$\ast\ast\ast\,\ast$} \label{eq:hard_lem_4}
    \end{align*}

    For the second part, we see that 
    \begin{align*}
        & \sum_{ab, uv \in E(G_{ij})} \sum_{k \neq i, j} k_{r-3}^{i,j,k} \abs{|N_k^c(u)| + |N^c_k(v)| - |N^c_k(a)| - |N_k^c(b)|} \\
        \leq{}& n^{r-3} \sum_{ab, uv \in E(G_{ij})} \sum_{k \neq i, j}  (\abs{|N_k^c(u)|- \ve n}  + \abs{|N^c_k(v)| - \ve n } + \abs{\ve n - |N^c_k(a)|}  + \abs{\ve n- |N_k^c(b)|}) \\
        \leq{}&  2n^r \sum_{k \neq i,j } \paren{\sum_{u \in V_i} \paren{\ve n - \abs{ N_k^c (u)}} + \sum_{v \in V_j} \paren{\ve n - \abs{ N_k^c (v)}}} \tag{$\ve n \geq \abs{N_k^c (u)}, \abs{N_k^c (v)}$} \\
        \leq{}&  4n^r \sum_{k \neq i,j } \paren{\ve n^2  - \paren{n^2 - \widehat{e(G)}}} \\
        \leq{}& 4  \ve^2 (r-1) (r-2) n^{r+2} \tag{$\ast\ast\ast\ast\ast$} \label{eq:hard_lem_5}
    \end{align*}
    where we used \Cref{eq:global_edge} in the last inequality. 

    Plugging the estimates from \eqref{eq:hard_lem_4} and \eqref{eq:hard_lem_5} into \eqref{eq:hard_lem_3} and summing over $ab, uv \in E(G_{ij})$, we derive 
    \begin{align*}
        & \sum_{ij \in \paren{[r] \setminus \set{1,2}}^{(2)}} \sum_{ab, uv\in E(G_{ij})} \abs{k_r(uv) - k_r(ab)} \\
        \leq{}& \binom{r-2}{2} \paren{2 \ve^2 (r-2) n^{r+2} + 4 \ve^2 (r-1) (r-2)n^{r+2} + 4 \ve^2 (r-2)^2 n^{r+2}}  \\
        \leq{}& 2(r-2) \paren{1 + 2(r-1) + 2(r-2) }\binom{r-2}{2} \ve^2 n^{r+2} \\
        \leq{}& 4 (r-2)^4 \ve^2 n^{r+2}. 
    \end{align*}
    Plugging this into \eqref{eq:hard_lem_2}, we conclude with (\ref{eq:hard_lem_1}) that
    \begin{align*}
        & \frac{1}{k_r(xy)} \sum_{R \in \K_r(xy)} {\delta_0(R- x-y ) }\\
        \leq{}& \frac{1}{(1- \const^{-1}/2) n^2} \paren{ \frac{3}{(1 - \const^{-1}) \widehat{e(G)}} \sum_{ij \in \paren{[r] \setminus \set{1,2}}^{(2)}} \sum_{ab, uv \in E(G_{ij})}  \abs{k_r(uv) - k_{r}(ab)} + \frac{2 \ve^2 r (r-2)^3 n^{r} }{ (1- \const^{-1})^2}} \\
        \leq{}&  \frac{1}{(1- \const^{-1}/2) n^2} \paren{ \frac{12 (r-2)^4 \ve^2 n^{r+2}}{(1 - \const^{-1}) \widehat{e(G)}} + \frac{2  \ve^2 r (r-2)^3 n^{r} }{ (1- \const^{-1})^2}} \\
        \leq{}& \frac{1}{(1- \const^{-1}/2) (1 - \const^{-1})}\paren{\frac{3n^2}{\widehat{e(G)}} + \frac{1}{2(1 - \const^{-1})}}4 \ve^2 (r-1)^2 (r-2)^2 n^{r-2} \\
        \leq{}& 15 \ve^2 (r-1)^2 (r-2)^2 n^{r-2} \tag{$\const \geq 26$}
    \end{align*}
    as claimed. 
\end{proof}
\subsection{Bounding edge-discrepancies}

\label{sec:disc}
We will now show that for $\omega_i$, the edge-discrepancy is decreasing at the desired rate. We will need the following lemma that, as with $h(\cdot, \cdot)$, one can show that the size of $\C(\cdot)$ is large and does not change much if varied among edges between two fixed parts.
\begin{lemma}
    For $e \in E(G)$, we have $\abs{\C(e)}\geq (1 - 3\ve)n^2$ for $n$ sufficiently large.
    Furthermore, for $e_{ij}^\ast, e_{ij} \in E(G_{ij})$, $1 \leq i < j \leq r$, $\abs{\abs{\C (e_{ij})} - \abs{\smash{\C(e_{ij}^\ast)}}} \leq 2 \varepsilon n^2$. 
    \label{clm:C-diff}
\end{lemma}
The proof is similar to \Cref{clm:e-diff} and \Cref{clm:h-diff} but we will include it here for completeness.

\begin{proof}
    For the first part, we directly compute 
    \begin{align*}
        \abs{\C(e)} \geq ((1- \ve) n-1) ((1 - 2 \ve ) n- 1) \overset{n \gg 1}{\geq } (1 - 3 \ve)n ^2.
    \end{align*}
    
    Without loss of generality $\abs{\C(e_{ij})} \geq \abs{\C(\smash{e_{ij}^\ast})}$. Let $u_i, u_i^\ast \in V_i, u_j, u_j^\ast \in V_j$ be the vertices of $e_{ij}, e_{ij}^\ast$. We see that 
    \begin{align*}
        & \abs{\C(e_{ij})} -  \abs{\C(e_{ij}^\ast)} \\
        ={}& \sum_{v_i \in V_i} \sum_{v_j \in V_j} \paren{\indicate{u_i u_j v_i v_j \in \C(e_{ij})} - \indicate{u_i^\ast u_j^\ast v_i v_j \in \C(e_{ij}^\ast)}} \\
        \leq{}& \sum_{v_i \in V_i \setminus\{u_i, u_i^*\}} \sum_{v_j \in V_j} \paren{\indicate{u_i u_j v_i v_j \in \C(e_{ij})} - \indicate{u_i^\ast u_j^\ast v_i v_j \in \C(e_{ij}^\ast)}} + 2\ve n \\
        \leq{}& \sum_{v_i \in V_i \setminus\{u_i, u_i^*\}} \sum_{v_j \in V_j \setminus \{u_j, u_j^*\}} \paren{\indicate{u_i u_j v_i v_j \in \C(e_{ij})} - \indicate{u_i^\ast u_j^\ast v_i v_j \in \C(e_{ij}^\ast)}} + 4\ve n \\
        \leq{}& \sum_{v_i \in N_i(u_j, u_j^\ast) \setminus \set{\smash{u_i^*, u_i}}} \sum_{v_j \in V_j \setminus \{u_j, u_j^*\}} \paren{\indicate{u_i u_j v_i v_j \in \C(e_{ij})} - \indicate{u_i^\ast u_j^\ast v_i v_j \in \C(e_{ij}^\ast)}} + \varepsilon n(n -2) + 4 \ve n\\
        \leq{}& \sum_{v_i \in N_i(u_j, u_j^\ast) \setminus \set{\smash{u_i^*, u_i}}} \sum_{v_j \in N_j(u_i, u_i^\ast) \setminus \set{ \smash{u_j^*, u_j}}} \paren{\indicate{u_i u_j v_i v_j \in \C(e_{ij})} - \indicate{u_i^\ast u_j^\ast v_i v_j \in \C(e_{ij}^\ast)}} + 2 \varepsilon n^2  \\ 
        ={}& 2\varepsilon n^2. 
    \end{align*}
    In the first inequality, note that the first indicator activates if $v_i = u_i^\ast$ at most $n-1$ times, while the second indicator activates if $v_i = u_i$ at least $(1- 2 \ve)n-1$ times; the reasoning for the second inequality is similar. The last equality follows because the vertices $v_i, v_j$ we sum over can only form $C_4$'s with $e_{ij}, e_{ij}^\ast$ respectively, if and only if $v_iv_j \in E(G)$.
\end{proof}

Another ingredient will be the following:
\begin{fact}\label{eq:sum_disc_vertex_zero}
    For distinct $i,j \in [r]$ and $v \in V_i$, we have for a vertex-balanced function $\omega_\ell$ that 
    \begin{align*}
        \sum_{u \in N_j(v)} \delta_{\ell} (uv) = 0. 
    \end{align*}
\end{fact}

\begin{proof}
    As $\omega_\ell$ is vertex-balanced by \Cref{cor:seq:v-balanced}, the claim follows from \Cref{rem:v-disc_e-disc}.
\end{proof}
\begin{lemma}\label{lem:conv}
    $\delta_{i+1}^\infty \leq 12\ve \delta_i^\infty$ provided that $\ve \leq 1/36$. In particular, $\lim_{i \to \infty} \delta_{i}^\infty = 0$ and $\omega_\infty$ is well-defined if $\ve \leq 1/36$. 
\end{lemma}
\begin{proof}
Let $e = v_kv_\ell$ be an edge lying between $V_k, V_\ell$. Fix $e_{k\ell}^\ast \in E(G_{k \ell})$, $v_k^\ast \in N_k(v_\ell)  \setminus \set{v_k}$, and $v_\ell^\ast \in N_\ell(v_k) \setminus \set{v_\ell}$. Using \Cref{eq:sum_disc_vertex_zero}, we obtain 
 \begin{align*}
     \delta_{i + 1}(e) ={}& \omega_{i + 1}(e) - \frac{k_r(G)}{\widehat{e(G)}}\\
    ={}& \omega_i(e) - \sum_{f \in E(G)} \delta_{i} (f) \psi_f(e) - \frac{k_r(G)}{\widehat{e(G)}}\\
    ={}& \delta_{i}(e) - \left(\sum_{f \in E(G)} \frac{\delta_{i} (f)}{\abs{\C(f) }} \sum_{C \in \C(f)}{\psi_{C , f}(e)}\right) \\
    ={}&  \delta_i(e) -  \sum_{f \in E(G_{k \ell} )} \frac{ (-1)^{\sigma(e,f)} \delta_i(f) |\C(e, f)|} {|\C(f)|} \tag{\Cref{fact:sum_C_f}}\\
    ={}& \delta_i(e) - \sum_{f \in E(G_{k \ell} ), |e \cap f| =0, 2}
    \frac{ \delta_i(f) |\C(e, f)|} {|\C(f)|} + \sum_{f \in E(G_{k \ell} ), |e \cap f| = 1}\frac{ \delta_i(f) |\C(e, f)|} {|\C(f)|}\\
    ={}& - \sum_{f \in E(G_{k \ell}), G[e\cup f] \simeq C_4} \frac{\delta_i(f)}{\abs{\C(f)}} + \sum_{f \in E(G_{k \ell} ), |e \cap f| = 1}\frac{ \delta_i(f) |\C(e, f)|} {|\C(f)|}. \tag{$\star$} \label{eq:conv_1} 
    \end{align*}
    
    Let us deal with the sums in \eqref{eq:conv_1} one by one: By \Cref{eq:sum_disc_vertex_zero}, $\sum_{e \in E(G_{k \ell})} \delta_i(e) = 0$, so
\begin{align*}
    & - \sum_{f \in E(G_{k \ell}), G[e\cup f] \simeq C_4} \frac{\delta_i(f)}{\abs{\C(f)}} \\
    ={}& - \sum_{f \in E(G_{k \ell}), G[e\cup f] \simeq C_4} \frac{\delta_i(f)}{\abs{\C(f)}} - \sum_{f \in E(G_{k \ell}),G[e\cup f] \not \simeq C_4 } \frac{\delta_i(f)}{|\C(e_{k\ell}^*)|} +\sum_{f \in E(G_{k \ell}), G[e\cup f] \not \simeq C_4 } \frac{\delta_i(f)}{|\C(e_{k\ell}^*)|} \\
    ={}&  \sum_{f \in E(G_{k \ell}), G[e\cup f] \simeq C_4} \paren{\frac{1}{\abs{\C(e_{k \ell}^\ast)}} - \frac{1}{\abs{\C(f)}}}\delta_i(f) + \sum_{f \in E(G_{k \ell}), G[e\cup f] \not \simeq C_4 } \frac{\delta_i(f)}{|\C(e_{k\ell}^*)|}.
\end{align*}
    Note that the number of $f \in E(G_{k\ell})$ such that $G[e \cup f] \not \simeq C_4$ is at most $2n-1 + 2 \ve n (n-1) = 2 \ve n^2 + \Oc(n)$. Similarly,
    we derive using \Cref{eq:sum_disc_vertex_zero}
\begin{align*}
    & {\sum_{f \in E(G_{k \ell} ), |e \cap f| = 1}\frac{ \delta_i(f) |\C(e, f)|} {|\C(f)|}} \\
    ={}& \sum_{u_k \in N_k(v_\ell) \setminus \set{v_k}}\frac{ \delta_i(v_\ell u_k) |\C(u_kv_\ell v_k )|} {|\C(u_kv_\ell)|} + \sum_{u_\ell \in N_\ell(v_k) \setminus \set{v_\ell}}\frac{ \delta_i(v_k u_\ell) |\C(u_\ell v_kv_\ell )|} {|\C(u_\ell v_k)|} \\
    ={}& \sum_{u_k \in N_k(v_\ell) \setminus \set{v_k}} \delta_i(v_\ell u_k) \left( \frac{|\C(u_kv_\ell v_k)|} {|\C(u_kv_\ell)|} -  \frac{|\C(v_k^*v_\ell v_k)|} {|\C(v_k^*v_\ell)|}\right) - \delta_{i}(v_k v_{\ell} ) \frac{|\C(v_k^*v_\ell v_k)|} {|\C(v_k^*v_\ell)|}\\
    & + \sum_{u_\ell \in N_\ell(v_k) \setminus \set{v_\ell}} \delta_i(v_k u_\ell) \left( \frac{  |\C(u_\ell v_k v_\ell )|} {|\C(u_\ell v_k)|}  - \frac{  |\C(v_\ell^* v_k v_\ell)|} {|\C(v_\ell^* v_k)|}\right) - \delta_{i}(v_k v_{\ell} ) \frac{|\C(v_\ell^*v_k v_\ell)|} {|\C(v_\ell^*v_k)|}.
\end{align*}

    Note that if $xy, yz \in E(G_{k \ell})$, where $x, z \in V_k$ are distinct, then 
    \begin{align*}
        \frac{\abs{\C(xyz)}}{\abs{\C(xy)}} &= \frac{\abs{N_\ell(x)} - \abs{N_\ell(x) \cap N_\ell^{c}(z)} - 1}{\sum_{a \in N_\ell(x) \setminus \set{y}} \paren{\abs{N_{k}(a,y)} - 1}} \\
        &\in \ofS{\frac{\abs{N_\ell(x)} -1 - \ve n }{\paren{\abs{N_\ell(x)} -1 }(n-1)}, \frac{\abs{N_\ell(x)} -1 }{\paren{\abs{N_\ell (x)}-  1}  ((1- 2\ve )n -1)}} \\
        &\subseteq \ofS{\frac{(1 -2 \ve )n - 1}{((1 - \ve )n - 1)(n-1)}, \frac{1}{(1- 2 \ve )n - 1} } \\
        &= \ofS{\frac{1-2\ve}{(1-\ve) n} + \Oc\paren{\frac{1}{n^2}}, \frac{1}{(1-2 \ve)n } + \Oc\paren{\frac{1}{n^2}}}
    \end{align*}
    
    Plugging this into \eqref{eq:conv_1}, we conclude with \Cref{clm:C-diff} that
    \begin{align*}
        & \abs{\delta_{i+1}(e)} \\
        \leq{}& \delta_i^\infty \paren{n^2 \paren{\frac{1}{(1 - 3\ve)n^2} - \frac{1}{n^2}} + \frac{2\ve n^2 + \Oc(n)}{(1- 3 \ve )n^2} + 2 \paren{n \ofS{\frac{1}{(1-2\ve)n} - \frac{1-2\ve }{(1-\ve)n}} + \Oc\paren{\frac{1}{n}}}} \\
        ={}& \delta_i^\infty  \paren{\frac{2 \varepsilon}{1 - \varepsilon} + \frac{2}{1 - 2 \varepsilon} + \frac{1 + 2 \varepsilon}{1 - 3 \varepsilon}  -  3 + \Oc\paren{\frac
         1 n }} \\
        ={}& \delta_i^{\infty}\left(\frac{2 \ve}{1 - \ve} + \frac{4 \ve}{1 - 2 \ve} + \frac{5 \ve}{1 - 3 \ve} + \Oc\paren{\frac
         1 n } \right)\\
        \leq{}& \delta_i^{\infty} \frac{11 \ve}{1 - 3 \ve}\\
        \leq{}& 12 \ve \delta_i^\infty 
    \end{align*}
    for $n$ sufficiently large and $\ve \leq 1/36$.

    That $\omega_\infty$ is well-defined follows then from \Cref{lem:omega_i+1_simplified}:
    \begin{align*}
        \|\omega_{i+1} - \omega_i\|_{\infty} &\leq \binom{r}{2} \delta_i^\infty \Oc\paren{n^{2-r}}
        \leq \binom{r}{2} \delta_0^\infty (12 \ve)^i \Oc\paren{n^{2-r}} .
    \end{align*}
    Thus, the difference geometrically tends to zero as $i$ tends to infinity.
\end{proof}

\subsection{Lower bounds on \texorpdfstring{$\omega_{i+1}(R)$}{weights}}

In this subsection, we will give lower bounds on $\omega_{i+1}(R)$, $R \in \K_r(G)$, compared to $\omega_{i}(R)$. This will be crucial in the proof of \Cref{thm:main}.

To have control on the size of $N(C)$ for a $4$-cycle $C \subseteq G_{ij}$, we first prove the following
\begin{lemma}
    If $C, C' \subseteq G_{ij}$ are $4$-cycles, then $\abs{\abs{N(C)} - \abs{N(C')}} \leq 4 \ve (r-2) n^{r-2}$. 
    \label{fact:N(C)_count}
\end{lemma}
\begin{proof}
    Without loss of generality $\abs{N(C)} \geq \abs{N(C')}$ and $i = 1, j=2$. 
    Then 
    \begin{align*}
        & \abs{\abs{N(C)} - \abs{N(C')}} \\
        \leq{}& \sum_{v_3 \in N_{3}(V(C))} \dots \sum_{v_r \in N_r(V(C))} \bigg(\indicate{v_3\dots v_r \in \K_{r-2}^{1,2} \text{ useful}} \\
        & - \indicate{v_3 \dots v_r \in \K_{r-2}^{1,2} \text{ useful}: v_3, \dots, v_r \in N(V(C'))}\bigg) \\
        \vdotswithin{=}& \\ 
        \leq{}& \sum_{v_3 \in N_{3}(V(C) \cup V(C'))} \dots \sum_{v_r \in N_r(V(C) \cup V(C'))} \bigg(\indicate{v_3\dots v_r \in \K_{r-2}^{1,2} \text{ useful}} \\
        & - \indicate{v_3 \dots v_r \in \K_{r-2}^{1,2} \text{ useful}:  v_3, \dots, v_r \in N(V(C'))}\bigg) + 4 \ve (r-2) n^{r-2} \\
        ={}& 4 \ve (r-2) n^{r-2}\qedhere
    \end{align*}
\end{proof}

The following estimate will be used for $i \geq 1$.
\begin{lemma}\label{lem:general_estimate_weighting}
    For all $i \geq 0$ and $R \in \K_r(G)$, $$\omega_{i+1} (R) \geq \omega_i(R) -  \frac{2r (r-1)}{\paren{1 - \frac{1}{\const } - \frac{2}{\constt}}\paren{1 - 3 \varepsilon}} \delta_i^\infty n^{2 - r}.$$
\end{lemma}
\begin{proof}
    By \Cref{lem:omega_i+1_simplified},
\begin{align*}
    \omega_{i+ 1}(R) 
    ={}& \omega_i(R) - \sum_{f \in E(R)}\sum_{C \in \C(f)}\sum_{e \in E(C)} 
       \frac{(-1)^{\sigma(e,f)}\delta_i(e)\indicate{\exists W \in N(C): G[W \cup f] = R}}{|N(C)||\C(e)|}\\
    \geq{}& \omega_i(R) -  \frac{1}{1 - \frac{1}{\const } - \frac{2}{\constt}} \sum_{f \in E(R)} \sum_{C \in \C(f)}\sum_{e \in E(C)} \frac{|\delta_i(e) |}{\paren{1 - 3 \varepsilon} n^r} \tag{\Cref{clm:C-diff}, \Cref{cor:useful_W}}\\
    \geq{}& \omega_i(R) - \frac{4}{1 - \frac{1}{\const } - \frac{2}{\constt}}  \sum_{f \in E(R)} \sum_{C \in \C(f)}\frac{\delta_i^\infty }{{1 - 3 \varepsilon}} \cdot n^{- r}\\
    \geq{}& \omega_i(R) -  \frac{2r (r-1)}{\paren{1 - \frac{1}{\const } - \frac{2}{\constt}}\paren{1 - 3 \varepsilon}} \delta_i^\infty n^{2 - r}. \qedhere{}
\end{align*}
\end{proof}

For $\omega_0$ to $\omega_1$, we will use the following 
\begin{lemma}\label{lem:omega_1_positive}
    For all $R \in \K_r(G)$ which are useful for some $f \in E(R)$, 
    \begin{equation*}
        \omega_1(R) \geq \omega_0(R) - \frac{1}{\paren{1 - \frac{1}{\const} - \frac{2}{\constt}}n^{r-2}} \abs{\delta_0 (R)}  - \binom{r}{2}  \frac{ 9(2r-3) \ve \delta_0^\infty n^{2-r}}{\paren{1- 3\ve}^2\paren{1 - \frac{1}{\const} - \frac{2}{\constt}}^2 }.
    \end{equation*}
\end{lemma}
\begin{proof}
    By \Cref{lem:omega_i+1_simplified},
\begin{align*}
        \omega_{1}(R) 
       &= \omega_0(R) - \sum_{f \in E(R)}\sum_{e \in E(G)}\sum_{C \in \C(f, e)}\delta_0(e)\frac{(-1)^{\sigma(e, f)}\indicate{\exists W \in N(C) : R = G[W \cup f]}}{|N(C)|{\abs{\C(e)}} }.
\end{align*}

For fixed $f = xy \in E(R)$ such that $R \in \widetilde{\K_r(f)}$, we break the sum over $e \in E(G)$ into four parts, $f = e$, $f \cap e = \{x\}$, $f \cap e = \{y\}$, and $f \cap e = \emptyset$. For all of them, we will assume that $x \in V_1$, $y \in V_2$ and that $R - f$ remains useful; otherwise, $f$ does not contribute to the sum. 

In the first case, we have using \Cref{cor:useful_W}
\begin{align*}
    &\abs{\sum_{e \in E(G), e = f}\sum_{C \in \C(f, e)}\delta_0(e)\frac{(-1)^{\sigma(e, f)}\indicate{\exists W \in N(C) : R = G[W \cup f]}}{{\abs{\C(e)}} |N(C)|}}\\ ={}& \abs{\sum_{C \in \C(f)}\delta_0(f)\frac{\indicate{\exists W \in N(C) : R = G[W \cup f]}}{{\abs{\C(f)}} |N(C)|}}\\
    \leq{}& \frac{1}{\paren{1- \frac{1}{\const} - \frac{2}{\constt}} n^{r - 2}} \abs{\delta_0(f)}.
\end{align*}

For the second case, fix an arbitrary edge $e^* \in E(G_{12}), e^* \neq f$ containing $x$, and some four cycle $C^*$ containing $x$ in $G_{12}$. 
\begin{align*}
    & \abs{\sum_{e \in E(G), e \cap f = \{x\}}\sum_{C \in \C(f, e)}\delta_0(e)\frac{(-1)^{\sigma(e, f)}\indicate{\exists W \in N(C) : R = G[W \cup f]}}{{\abs{\C(e)}} |N(C)|}}\\ 
    ={}& \abs{\sum_{v \in N_2(x) \setminus\{y\}}\sum_{C \in \C(v, f)}\delta_0(vx)\frac{\indicate{\exists W \in N(C) : R = G[W \cup f]}}{{\abs{\C(vx)}} |N(C)|}}\\
    \leq{}& \abs{\sum_{v \in N_2(x) \setminus\{y\}}\sum_{C \in \C(v, f)}\delta_0(vx)\frac{\indicate{\exists W \in N(C) : R = G[W \cup f]}}{{\abs{\C(vx)}} |N(C^*)|}} \\
    &+ \abs{\sum_{v \in N_2(x) \setminus\{y\}}\sum_{C \in \C(v, f)}\delta_0(vx)\frac{\indicate{\exists W \in N(C) : R = G[W \cup f]}\paren{|N(C)| - |N(C^*)|}}{{\abs{\C(vx)}} |N(C)||N(C^*)|}}  \\ 
    \leq{}& \abs{\sum_{v \in N_2(x) \setminus\{y\}}\sum_{C \in \C(v, f)}\delta_0(vx)\frac{\indicate{\exists W \in N(C) : R = G[W \cup f]}}{{\abs{\C(vx)}} |N(C^*)|}} + \frac{4 \ve (r-2) \delta_0^{\infty} n^{2 - r}}{\paren{1 - 3 \ve }\paren{1 - \frac{1 }{\const} - \frac{2}{\constt}}^2} \tag{\Cref{cor:useful_W}, \Cref{clm:C-diff}, \Cref{fact:N(C)_count}}\\
    \leq{}& \abs{\sum_{v \in N_2(x) \setminus\{y\}}\sum_{C \in \C(v, f)}\delta_0(vx)\frac{\indicate{\exists W \in N(C) : R = G[W \cup f]}}{{\abs{\C(e^*)}} |N(C^*)|}}  \\&+ \abs{\sum_{v \in N_2(x) \setminus\{y\}}\sum_{C \in \C(v, f)}\delta_0(vx)\frac{\indicate{\exists W \in N(C): R = G[W \cup f]} (|\C(vx)| - |\C(e^*)|)}{\abs{\C(vx)}\abs{\C(e^*)} |N(C^*)|}} \\
    & +  \frac{4 \ve (r-2) \delta_0^{\infty} n^{2 - r}}{\paren{1 - 3 \ve }\paren{1 - \frac{1 }{\const} - \frac{2}{\constt}}^2} \\
    \leq{}& \abs{\sum_{v \in N_2(x) \setminus\{y\}}\sum_{C \in \C(v, f)}\delta_0(vx)\frac{\indicate{\exists W \in N(C) :  R = G[W \cup f]}}{\abs{\C(e^*)}|N(C^*)|}} + \frac{(4r-6) \ve \delta_0^{\infty} n^{2 - r}}{\paren{1 - 3 \ve }^2\paren{1 - \frac{1 }{\const} - \frac{2}{\constt}}^2} \tag{\Cref{cor:useful_W}, \Cref{clm:C-diff}}\\
    \leq{}& \abs{\sum_{v \in N_2(x) \setminus\{y\}}\sum_{C \in \C(v, f)}\frac{\delta_0(vx)}{\abs{\C(e^*)}|N(C^*)|}}+ {\sum_{v \in N_2(x) \setminus\{y\}}\sum_{C \in \C(v, f)}\abs{\delta_0(vx)}\frac{\indicate{V(C)\not \subseteq N(V(R) \setminus f)}}{\abs{\C(e^*)}|N(C^*)|}}  \\& + \frac{(4r-6) \ve \delta_0^{\infty} n^{2 - r}}{\paren{1 - 3 \ve }^2\paren{1 - \frac{1 }{\const} - \frac{2}{\constt}}^2} \\
     \leq{}& \abs{\sum_{v \in N_2(x) \setminus\{y\}}\sum_{C \in \C(v, f)}\frac{\delta_0(vx)}{\abs{\C(e^*)}|N(C^*)|}} + \frac{2 (r-2) \ve n^2 \delta_0^\infty}{\abs{\C(e^\ast)} \abs{N(C^\ast)}}  + \frac{(4r-6) \ve \delta_0^{\infty} n^{2 - r}}{\paren{1 - 3 \ve }^2\paren{1 - \frac{1 }{\const} - \frac{2}{\constt}}^2} \\
     \leq{}& \abs{\sum_{v \in N_2(x) \setminus\{y\}}\sum_{C \in \C(v, f)}\frac{\delta_0(vx)}{\abs{\C(e^*)}|N(C^*)|}}  + \frac{(6r-10) \ve \delta_0^{\infty} n^{2 - r}}{\paren{1 - 3 \ve }^2\paren{1 - \frac{1 }{\const} - \frac{2}{\constt}}^2}  \tag{\Cref{cor:useful_W}}\\
    ={}& \abs{\sum_{v \in N_2(x) \setminus\{y\}}\frac{|\C(vxy)|\delta_0(vx)}{\abs{\C(e^*)}|N(C^*)|}}  + \frac{(6r-10) \ve \delta_0^{\infty} n^{2 - r}}{\paren{1 - 3 \ve }^2\paren{1 - \frac{1 }{\const} - \frac{2}{\constt}}^2}  \\
    \leq{}& \abs{\sum_{v \in N_2(x) \setminus\{y\}}\frac{|\C(e^\ast y)|\delta_0(vx)}{\abs{\C(e^*)}|N(C^*)|}}  + \abs{\sum_{v \in N_2(x) \setminus\{y\}}\frac{(|\C(vxy)| - \abs{\C(e^\ast y )})\delta_0(vx)}{\abs{\C(e^*)}|N(C^*)|}}  + \frac{(6r-10) \ve \delta_0^{\infty} n^{2 - r}}{\paren{1 - 3 \ve }^2\paren{1 - \frac{1 }{\const} - \frac{2}{\constt}}^2}  \\
    \leq{}& \abs{\sum_{v \in N_2(x) \setminus\{y\}}\frac{|\C(e^*y)|\delta_0(vx)}{\abs{\C(e^*)}|N(C^*)|}}  + \frac{(6r-10) \ve \delta_0^{\infty} n^{2 - r}}{\paren{1 - 3 \ve }^2\paren{1 - \frac{1 }{\const} - \frac{2}{\constt}}^2} + \frac{\ve \delta_0^\infty n^{2-r}}{(1- 3\ve)\paren{1 - \frac{1}{\const} - \frac{2}{\constt}}}\tag{\Cref{cor:useful_W}, \Cref{clm:C-diff}}  \\
    \leq{}& \abs{\sum_{v \in N_2(x)}\frac{|\C(e^*y)|\delta_0(vx)}{\abs{\C(e^*)}|N(C^*)|}}  + \frac{3 (2r-3) \ve \delta_0^{\infty} n^{2 - r}}{\paren{1 - 3 \ve }^2\paren{1 - \frac{1 }{\const} - \frac{2}{\constt}}^2} \tag{$n \gg 1$}\\
    \leq{}&   \frac{3 (2r-3) \ve \delta_0^{\infty} n^{2 - r}}{\paren{1 - 3 \ve }^2\paren{1 - \frac{1 }{\const} - \frac{2}{\constt}}^2}. \tag{\Cref{eq:sum_disc_vertex_zero}}
\end{align*}
The third case follows symmetrically. 

For the fourth case, choose $e^\ast \in E(G_{12})$ with $e^\ast \cap f = \emptyset$ and $C^\ast \in \C(f)$: 
\begin{align*}
    & \abs{\sum_{e \in E(G_{12}) \colon e \cap f = \emptyset }\sum_{C \in \C(f, e)}\delta_0(e)\frac{(-1)^{\sigma(e, f)}\indicate{\exists W \in N(C) : R = G[W \cup f]}}{|N(C)|{\abs{\C(e)}} }}\\
    ={}& \abs{\sum_{x' y' \in E(G_{12}), G[xy x'y'] \simeq C_4} \delta_0(x'y')\frac{\indicate{x',y' \in N(V(R) \setminus f)}}{|N(xyx'y')|{\abs{\C(x'y')}} } } \\
    \leq{}& \abs{\sum_{x' y' \in E(G_{12}), G[xy x'y'] \simeq C_4} \delta_0(x'y')\frac{\indicate{x',y' \in N(V(R) \setminus f)}}{|N(C^\ast)|{\abs{\C(x'y')}} } } \\
    & + \abs{\sum_{x' y' \in E(G_{12}), G[xy x'y'] \simeq C_4} \delta_0(x'y')\frac{\indicate{x',y' \in N(V(R) \setminus f)} \paren{\abs{N(C^\ast)} - \abs{N(xyx'y')}}}{\abs{N(xyx'y')}|N(C^\ast)|{\abs{\C(x'y')}} } } \\
    \leq{}& \abs{\sum_{x' y' \in E(G_{12}), G[xy x'y'] \simeq C_4} \delta_0(x'y')\frac{\indicate{x',y' \in N(V(R) \setminus f)}}{|N(C^\ast)|{\abs{\C(x'y')}} } } + \frac{4 \ve (r-2) \delta_0^\infty n^{2-r}}{\paren{1- 3\ve}\paren{1 - \frac{1}{\const} - \frac{2}{\constt}}^2 } \tag{\Cref{cor:useful_W}, \Cref{clm:C-diff}, \Cref{fact:N(C)_count}}\\
    \leq{}& \abs{\sum_{x' y' \in E(G_{12}), G[xy x'y'] \simeq C_4} \delta_0(x'y')\frac{\indicate{x',y' \in N(V(R) \setminus f)}}{|N(C^\ast)|{\abs{\C(e^\ast)}} } } + \frac{4 \ve (r-2) \delta_0^\infty n^{2-r}}{\paren{1- 3\ve}\paren{1 - \frac{1}{\const} - \frac{2}{\constt}}^2 }  \\
    & + \abs{\sum_{x' y' \in E(G_{12}), G[xy x'y'] \simeq C_4} \delta_0(x'y')\frac{\indicate{x',y' \in N(V(R) \setminus f)}(\abs{\C(e^\ast)} - \abs{\C(x'y')})}{|N(C^\ast)|{\abs{\C(e^\ast)}} \abs{\C (x'y')} } } \\
    \leq{}& \abs{\sum_{x' y' \in E(G_{12}), G[xy x'y'] \simeq C_4} \delta_0(x'y')\frac{\indicate{x',y' \in N(V(R) \setminus f)}}{|N(C^\ast)|{\abs{\C(e^\ast)}} } } + \frac{ (4r-6) \ve \delta_0^\infty n^{2-r}}{\paren{1- 3\ve}^2\paren{1 - \frac{1}{\const} - \frac{2}{\constt}}^2 } \tag{\Cref{cor:useful_W}, \Cref{clm:C-diff}} \\
    \leq{}& \abs{\sum_{x' y' \in E(G_{12}), G[xy x'y'] \simeq C_4} \frac{\delta_0(x'y')}{|N(C^\ast)|{\abs{\C(e^\ast)}} } }  + \frac{ (4r-6) \ve \delta_0^\infty n^{2-r}}{\paren{1- 3\ve}^2\paren{1 - \frac{1}{\const} - \frac{2}{\constt}}^2 } \\
    & + \abs{\sum_{x' y' \in E(G_{12}), G[xy x'y'] \simeq C_4} \delta_0(x'y')\frac{\indicate{\set{x',y'} \not \subseteq N(V(R) \setminus f)}}{|N(C^\ast)|{\abs{\C(e^\ast)}} } }  \\
    \leq{}& \abs{\sum_{x' y' \in E(G_{12}), G[xy x'y'] \simeq C_4} \frac{\delta_0(x'y')}{|N(C^\ast)|{\abs{\C(e^\ast)}} } }  + \frac{ (6r-10) \ve \delta_0^\infty n^{2-r}}{\paren{1- 3\ve}^2\paren{1 - \frac{1}{\const} - \frac{2}{\constt}}^2 } \tag{\Cref{cor:useful_W}, \Cref{clm:C-diff}}  \\
    \leq{}& \abs{\sum_{x' \in N_1(y) \setminus \set{x}} \sum_{y' \in N_2(x')} \frac{\delta_0(x'y')}{|N(C^\ast)|{\abs{\C(e^\ast)}} } } + \frac{ (6r-10) \ve \delta_0^\infty n^{2-r}}{\paren{1- 3\ve}^2\paren{1 - \frac{1}{\const} - \frac{2}{\constt}}^2 }\\
    & + \abs{\sum_{x' \in N_1(y) \setminus \set{x}} \sum_{y' \in N_2(x') \setminus \set{y}} \frac{\delta_0(x'y') \indicate{y' \not \in N_2(x)}}{|N(C^\ast)|{\abs{\C(e^\ast)}} } }  + \abs{\sum_{x' \in N_1(y) \setminus \set{x}} \frac{\delta_0(x'y)}{\abs{N(C^\ast) }\abs{\C(e^\ast)}}}\\
    \leq{}& \abs{\sum_{x' \in N_1(y) \setminus \set{x}} \sum_{y' \in N_2(x')} \frac{\delta_0(x'y')}{|N(C^\ast)|{\abs{\C(e^\ast)}} } } + \frac{ 3(2r-3) \ve \delta_0^\infty n^{2-r}}{\paren{1- 3\ve}^2\paren{1 - \frac{1}{\const} - \frac{2}{\constt}}^2 } \tag{\Cref{cor:useful_W}, \Cref{clm:C-diff}, $n \gg 1$}  \\
    ={}&  \frac{ 3(2r-3) \ve \delta_0^\infty n^{2-r}}{\paren{1- 3\ve}^2\paren{1 - \frac{1}{\const} - \frac{2}{\constt}}^2 }.\tag{\Cref{eq:sum_disc_vertex_zero}}
\end{align*}

Combining all of these bounds gives the desired conclusion.
\end{proof}

\section{Proof of \texorpdfstring{\Cref{thm:main}}{main theorem}}\label{sec:proof}

\begin{proof}[Proof of \Cref{thm:main}]
    Let $G$ be a balanced, $r$-partite, $K_r$-divisible graph with $\hat{\delta}(G) \geq (1- \ve)n$, $n \gg 1$. Without loss of generality we may assume that $G$ is $(\ve, r)$-minimal. It suffices to show that, after normalization, $\omega_{\infty}$ is a fractional $K_r$-decomposition. By \Cref{lem:conv}, the only thing remaining to check is that $\omega_{\infty}$ is non-negative. Let $R \in \K_r(G)$ be an arbitrary $r$-clique. 

    By \Cref{lem:omega_0_positive}, 
    \begin{equation*}
        \omega_0(R) \geq 1 - \frac{1}{\const - 1 }.
    \end{equation*}
    If $R$ is \emph{not} useful for any edge, then this shows that $\omega_{\infty} (R) = \omega_0(R) > 0$ since we only build edge-adjusters using useful $r$-cliques. So, assume now that $R$ is useful for some edge $e$, that is, 
    \begin{equation*}
        \delta_0(R) \leq  \paren{\frac{15 \constt}{\const} + 13}\ve (r-1)(r-2) n^{r-2}.
    \end{equation*}
    
    By \Cref{lem:omega_1_positive}, we get after the first edge adjustment step
    \begin{align*}
        \omega_1(R) &\geq 1 - \frac{1}{\const - 1 } - \frac{1}{\paren{1 - \frac{1}{\const} - \frac{2}{\constt}}n^{r-2}} \abs{\delta_0 (R)}  - \binom{r}{2}  \frac{ 9(2r-3) \ve \delta_0^\infty n^{2-r}}{\paren{1- 3\ve}^2\paren{1 - \frac{1}{\const} - \frac{2}{\constt}}^2 } \\
        &\geq 1 - \frac{1}{\const - 1 } - \paren{\frac{15 \constt}{\const} + 13}\frac{\ve (r-1)(r-2) n^{r-2} }{\paren{1 - \frac{1}{\const} - \frac{2}{\constt}}n^{r-2}}  - \binom{r}{2}  \frac{ 117(2r-3) (r-2) \ve^2 }{2\paren{1- 3\ve}^2\paren{1 - \frac{1}{\const} - \frac{2}{\constt}}^2 } \\
        &\geq 1 - \frac{1}{\const - 1 } - \paren{\frac{15 \constt}{\const} + 13}\frac{1 }{\paren{\const - 1- \frac{2\const }{\constt}}}  - \frac{117}{2\paren{1- 3\ve}^2\paren{\const  - 1 - \frac{2\const}{\constt}}^2 } 
    \end{align*}
    where we used \Cref{genericvertexdiscrepencybound} and $R$ being useful in the second inequality.

    Thus, we obtain for $\omega_{\infty}(R)$ using \Cref{genericvertexdiscrepencybound}, \Cref{lem:general_estimate_weighting}, and \Cref{lem:conv}  
\begin{align*}
    \omega_\infty(R) &\geq \omega_1(R)  - \frac{2r (r-1)}{\paren{1 - \frac{1}{\const } - \frac{2}{\constt}}\paren{1 - 3 \varepsilon}} n^{2 - r} \sum_{i = 1}^\infty \delta_i^\infty  \\ 
    &\geq \omega_1(R) -  \frac{13 r (r-1) (r-2) \ve}{\paren{1 - \frac{1}{\const } - \frac{2}{\constt}}\paren{1 - 3 \varepsilon}}  \sum_{i = 1}^\infty \paren{12 \ve}^i \tag{$\ve \leq 1/36$} \\
    &= \omega_1(R) -  \frac{13 r (r-1) (r-2) \ve}{\paren{1 - \frac{1}{\const } - \frac{2}{\constt}}\paren{1 - 3 \varepsilon}}  \frac{12 \ve}{1 - 12 \ve }  \\
    &\geq \omega_1(R) -  \frac{156 \ve (r-2) }{\paren{\const - 1 - \frac{2 \const}{\constt}}\paren{1 - 15 \varepsilon}}  \\
    &\geq 1 - \frac{1}{\const - 1 } - \frac{1}{\const - 1- \frac{2\const }{\constt}} \paren{\frac{15 \constt}{\const} + 13 + \frac{117}{2\paren{1- 3\ve}^2\paren{\const  - 1 - \frac{2\const}{\constt}} } +  \frac{156 \ve (r-2) }{{1 - 15 \varepsilon}}} \\
    &\geq 1 - \frac{1}{29} - \frac{1}{23} \paren{18 + \frac{117}{46(1- 3 \ve)^2 } + \frac{156 \ve (r-2)}{1- 15 \ve }} \\
    &\geq 0 
\end{align*}
for $r \geq 3$ provided that $\const = 30, \constt = 10$.   
\end{proof}

\section{Acknowledgements}
All statements and proofs were developed and completed by the authors in the traditional way, without any use of AI tools. The paper was written entirely by the authors. The authors used ChatGPT (GPT-6.1 Sol) for assistance with proofreading only.  It also pointed out that the constant $30$ in \Cref{thm:main} may be replaced by $18$ by being more efficient in some of our arguments but we decided to keep our original proof.

\bibliographystyle{amsplain_noetal}
\bibliography{refs}

\end{document}